\documentclass[11pt,oneside]{amsart}
\usepackage[top=25mm, bottom=25mm, left=20mm, right=20mm]{geometry}
\usepackage{amsfonts,amssymb,amsmath,amsthm,amsxtra}
\usepackage[T1]{fontenc}
\usepackage[utf8]{inputenc}
\usepackage{bm}
\usepackage{mathtools}
\usepackage{dsfont}
\usepackage{bbm}
\usepackage{mathrsfs}
\usepackage{enumitem}
\usepackage{color}

\usepackage{graphics}

\numberwithin{equation}{section}
\newtheorem{theorem}{Theorem}
\newtheorem{corollary}[theorem]{Corollary}
\newtheorem{lemma}[theorem]{Lemma}
\newtheorem{proposition}[theorem]{Proposition}
\theoremstyle{definition}

\newcommand{\BB}{\mathbb{B}}
\newcommand{\CC}{\mathbb{C}}

\newcommand{\GG}{\mathbb{G}}

\newcommand{\II}{\mathbb{I}}

\newcommand{\LL}{\mathbb{L}}

\newcommand{\NN}{\mathbb{N}}

\newcommand{\PP}{\mathbb{P}}
\newcommand{\QQ}{\mathbb{Q}}
\newcommand{\RR}{\mathbb{R}}

\newcommand{\TT}{\mathbb{T}}

\newcommand{\ZZ}{\mathbb{Z}}

\newcommand{\N}{\mathbb{N}}

\newcommand{\calA}{\mathcal{A}}

\newcommand{\calP}{\mathcal{P}}
\newcommand{\calF}{\mathcal{F}}
\newcommand{\calB}{\mathcal{B}}
\newcommand{\calM}{\mathcal{M}}

\newcommand{\calH}{\mathcal{H}}

\newcommand{\calK}{\mathcal{K}}

\newcommand{\calS}{\mathcal{S}}

\newcommand{\calQ}{\mathcal{Q}}

\newcommand{\frakm}{\mathfrak{m}}

\newcommand{\frkS}{\mathfrak{S}}

\newcommand{\frkM}{\mathfrak{m}}
\newcommand{\frkN}{\mathfrak{n}}
\newcommand{\frkY}{\mathfrak{y}}

\newcommand{\frkI}{\mathfrak{i}}

\newcommand{\dif}{\mathrm{d}}
\newcommand{\ex}{\bm{e}}

\newcommand{\ind}[1]{\mathds{1}_{{#1}}}

\newcommand*{\DMO}[1]{\expandafter\DeclareMathOperator\csname #1\endcsname {#1}}
\DMO{ord}
\DMO{supp}
\DMO{Sym}
\DMO{cl}
\DMO{lcm}
\DMO{dist}
\DMO{tr}
\DMO{diam}

\DeclarePairedDelimiter\abs{\lvert}{\rvert}

\DeclarePairedDelimiter\norm{\lVert}{\rVert}

\DeclarePairedDelimiterX\spr[2]{\langle}{\rangle}{#1,#2}
\newcommand{\ipr}[2]{#1\cdot#2}
\DeclarePairedDelimiterX\Set[2]{\{}{\}}{#1\colon #2}
\DeclarePairedDelimiterX\Seq[1]{(}{)}{#1}

\begin{document}
\title[Vector-valued maximal inequalities for polynomial ergodic averages]{Vector-valued maximal inequalities and multi-parameter oscillation inequalities for the polynomial ergodic averages along multi-dimensional subsets of primes}

\author{Nathan Mehlhop}

\address[Nathan Mehlhop]{Department of Mathematics,
Rutgers University,
Piscataway, NJ 08854-8019, USA }
\email{nam225@math.rutgers.edu}
\thanks{The author was partially supported by the National Science Foundation (NSF) grant DMS-2154712.}

\begin{abstract}
We prove the uniform $\ell^2$-valued maximal inequalities for polynomial ergodic averages and truncated singular operators of Cotlar type modeled over multi-dimensional subsets of primes. In the averages case, we combine this with earlier one-parameter oscillation estimates \cite{MeS} to prove corresponding multi-parameter oscillation estimates. This provides a fuller quantitative description of the pointwise convergence of the mentioned averages and is a generalization of the polynomial Dunford-Zygmund ergodic theorem attributed to Bourgain \cite{MISZWR}.
\end{abstract}
\maketitle
\section{Introduction}
\subsection{Statement of results}
 Let $(X,\calB,\mu)$ be a  $\sigma$-finite measure space endowed with a family of invertible commuting and measure preserving transformations $S_1,\ldots, S_d:X\to X$. Let $\Omega$ be a bounded convex open subset of $\RR^k$ such that $B(0,c_{\Omega}) \subseteq \Omega \subseteq B(0,1)$ for some $c_{\Omega}\in(0, 1)$, where $B(0, u)$ is the open Euclidean ball in $\RR^k$ with radius $u>0$ centered at $0\in\RR^k$. For any $t>0$, we set
\[
\Omega_t := \{x\in\RR^{k}: t^{-1}x\in\Omega\}.
\] 
We consider a polynomial mapping
\begin{equation}\label{polymap}
    \mathcal{P}=(\mathcal{P}_1,\dots,\mathcal{P}_{d})\colon\ZZ^k\to\ZZ^{d}
\end{equation}
where each $\calP_j\colon\ZZ^k\to\ZZ$ is a polynomial of $k$ variables with integer coefficients such that $\calP_j(0)=0$. Let $k',k''\in\{0,1,\ldots,k\}$ with $k=k'+k''$. For $f\in L^\infty(X,\mu)$, we define the associated ergodic averages by  
\begin{equation}\label{def:average}
\calA_{t}^{\calP,k',k''}f(x):=\frac{1}{\vartheta_\Omega(t)}\sum_{(n,p)\in\ZZ^{k'}\times(\pm\PP)^{k''}}f(S_1^{\calP_1(n,p)}\cdots S_d^{\calP_d(n,p)}x)\mathds{1}_{\Omega_t}(n,p)\Big(\prod_{i=1}^{k''}\log|p_i|\Big),\quad x\in X,
\end{equation}
where $\pm\PP$ denotes the set of positive and negative prime numbers and
\begin{equation*}
\vartheta_\Omega(t):=\sum_{(n,p)\in\ZZ^{k'}\times(\pm \PP)^{k''}}\mathds{1}_{\Omega_t}(n,p)\Big(\prod_{i=1}^{k''}\log|p_i|\Big)
\end{equation*}
is the Chebyshev function. We also consider the Cotlar type ergodic averages given by
\begin{equation}\label{def:singular}
\calH_{t}^{\calP,k',k''}f(x):=\sum_{(n,p)\in\ZZ^{k'}\times(\pm\PP)^{k''}}f(S_1^{\calP_1(n,p)}\cdots S_d^{\calP_d(n,p)}x)K(n,p)\mathds{1}_{\Omega_t}(n,p)\Big(\prod_{i=1}^{k''}\log|p_i|\Big),\quad x\in X,
\end{equation}
where $K\colon\RR^{k}\setminus\{0\} \to \CC$ is a Calder\'on--Zygmund kernel satisfying the following conditions:
\begin{enumerate}
\item The size condition: For every $x\in\RR^k\setminus\{0\}$, we have
\begin{equation}
\label{eq:size-unif}
\abs{K(x)} \lesssim \abs{x}^{-k}.
\end{equation}
\item The cancellation condition: For every $0<r<R<\infty$, we have
\begin{equation}
\label{eq:cancel}
\int_{\Omega_{R}\setminus \Omega_{r}} K(y) \dif y = 0.
\end{equation}
\item
  The Lipschitz continuity condition:
  For every $x, y\in\RR^k\setminus\{0\}$ with $2|y|\leq |x|$, we have
\begin{equation}
\label{eq:K-modulus-cont}
\abs{K(x)-K(x+y)}
\lesssim 
\abs{y} \abs{x}^{-(k+1)}.
\end{equation}
\end{enumerate}
For a sequence of functions $(f_\frkI)_{\frkI\in\NN}$ with each $f_\frkI \in L^p(X,\mu)$, we define the $L^p(X;\ell^2)$ norm by
\begin{equation}
    \|f_\frkI\|_{L^p(X;\ell^2)} = \Big\|\big(\sum_{\frkI \in \NN} |f_\frkI|^2\big)^{1/2}\Big\|_{L^p(X)}
\end{equation}
and we say that $(f_\frkI)_{\frkI\in\NN} \in L^p(X;\ell^2)$ if $\|f_\frkI\|_{L^p(X;\ell^2)} < \infty$.

We can now state the main result of this paper.
\begin{theorem}\label{MAINthm}
Let $d, k \in \NN$ and let $\calP$ be a polynomial mapping as in ~\eqref{polymap}. Let $k',k''\in\{0,1,\ldots,k\}$ with $k'+k''=k$ and let $\calM_t^{\calP,k',k''}$ be either $\calA_t^{\calP,k',k''}$ or $\calH_t^{\calP,k',k''}$. Then, for any $p\in(1,\infty)$, there is a constant $C_{p,d,k,\deg\calP}>0$ such that
\begin{equation}
\label{vvest}
    \big\|\sup_{t > 0}|\calM_t^{\calP,k',k''} f_\frkI|\big\|_{L^p(X;\ell^2)} \leq C_{p,d,k,\deg \calP}\|f_\frkI\|_{L^p(X;\ell^2)}
\end{equation}
for any $(f_\frkI)_{\frkI \in \NN} \in L^p(X;\ell^2)$. The constant $C_{p,d,k,\deg\calP}$ is independent of the
coefficients of the polynomial mapping $\calP$.
\end{theorem}

In the proof of the above theorem, we use methods developed in \cite{MST2,MSZ3,Troj} and very recently in \cite{MeS,MSS,S}. We follow Bourgain's approach \cite{B3} to use the Calderón transference principle \cite{Cald} which reduces the problem to the integer shift system (see Section~\ref{sec:red}) and then exploit the Hardy--Littlewood circle method to analyze the appropriate Fourier multipliers. The main tools used to handle the estimates for the multiplier operators are: an appropriate generalization of Weyl's inequality (Proposition~\ref{weylinq}); the Ionescu--Wainger multiplier theorem (see \cite{IW,MSZ3} and \cite{TaoIW}) combined with the Rademacher--Menshov inequality (see \cite{MST2}) and standard multiplier approximations (Lemma~\ref{multiplierapprox}); the Magyar--Stein--Wainger sampling principle \cite{MSW} and \cite{MSZ1}. Throughout, we also use the Marcinkiewicz-Zygmund inequality (Proposition \ref{MZprop}) to extend scalar inequalities to their vector-valued analogues. 

We recall the $\lambda$-jump
counting function and the variation and oscillation semi-norms, which give quantitative measures for pointwise convergence. We use the
convention that a supremum taken over the empty set is zero. Let $\II \subseteq \RR$ with $\#\II \geq 2$ and $f: \II \rightarrow \CC$. Given $\lambda > 0$,the \emph{$\lambda$-jump counting function} of $f$ is defined by
\begin{equation*}
  \label{def:jump}
N_{\lambda}(f(t) : t\in\II):=\sup \{J\in\NN\, |\, \exists_{\substack{t_{0}<\dotsb<t_{J}\\ t_{j}\in\II}}  : \min_{0<j\leq J} \abs{f(t_{j})-f(t_{j-1})} \geq \lambda\}.
\end{equation*}
Given $r\in[1,\infty)$, the $r$-\textit{variation seminorm} $V^r$ of a $f$ is defined by
\begin{equation*}\label{def:var}
	V^r(f(t) : t \in \II): = \sup_{\substack{t_0 < t_1 < \cdots < t_N \\ t_j \in \II}}
	\Big(\sum_{j = 1}^N \abs{f(t_j) - f(t_{j-1})}^r\Big)^{1/r}.
\end{equation*}
Given $r \in [1,\infty)$, $N\in\NN\cup\{\infty\}$, and a strictly increasing sequence $I=(I_j: j\in\NN)\subseteq\II$, the {\it truncated $r$-oscillation seminorm} of $f$ is defined by
\begin{equation*}
\label{def:osc}
O_{I, N}^r( f(t): t\in \II)
:= \Big(\sum_{j=1}^N\sup_{\substack{I_j \le t < I_{j+1}\\t\in\II}}
\abs{f(t)-f(I_j)}^r\Big)^{1/r}.
\end{equation*}
Because of their preeminent role in multi-parameter pointwise convergence problems, we also consider multi-parameter analogues of the oscillation seminorms. Let $\II \subseteq \RR^M$ with $\#\II \geq 2$ and $f: \II \rightarrow \CC$. Given $r \in [1,\infty)$, $N\in\NN\cup\{\infty\}$, and a sequence $I=(I_j: j\in\NN)\subseteq\II$ that is strictly increasing in every coordinate, the $M$-\textit{parameter truncated r-oscillation seminorm} of $f$ is defined by
\begin{equation*}
O_{I, N}^r(f(t): t \in \II):=
\Big(\sum_{j=1}^{N}\sup_{t\in \BB[I_j]\cap\II}
\abs{f(t) - f(I_j)}^r\Big)^{1/r},
\end{equation*}
where
$\BB[I_j]:=[I_{j1}, I_{(j+1)1})\times\ldots\times[I_{jM}, I_{(j+1)M})$
is a box determined by the element $I_j=(I_{j1}, \ldots, I_{jM})$ of
the sequence $I$. Without causing any ambiguity, we may instead take $I \in \frkS_N(\II)$, the family of all 
sequences of length $N+1$ contained in $\II$ that are strictly increasing in every coordinate. For more information about these quantitative tools in the study of pointwise convergence problems, we refer to \cite{MSWr}, see also \cite{B3,jsw,MSS,MISZWR,MSZ1,S1}.

We now recall an abstract multi-parameter oscillation result. For a linear operator $T: L^0(X) \rightarrow L^0(X)$, we denote by $|T|$ the sublinear maximal operator taken in the lattice sense defined by 
\[|T|f(x) = \sup_{|g| \leq |f|} |Tg(x)|, \quad x \in X, \quad f \in L^p(X).\]

\begin{proposition} \label{prop:multi} \cite[Proposition 4.1]{MSWr}
Let $(X,\mathcal B(X), \mu)$ be a $\sigma$-finite measure space and let $\II\subseteq \RR$ be such that $\#\II\ge2$.
Let $k \in \NN_{\ge2}$ and  $p, r\in(1, \infty)$ be fixed. Let $(T_t)_{t\in\II^k}$ be a family of linear operators of the form
\begin{align*}
T_t :=  T_{t_1}^{1} \cdots T_{t_k}^{k},
\qquad t=(t_1,\ldots, t_k)\in\II^k,
\end{align*}
where
$\{ T_{t_i}^i : i \in [k], \, t_i \in \II \}$ is a family of commuting linear operators that are bounded on $L^p(X)$.
If the set $\II$ is uncountable, then we also assume that $\II\ni t\mapsto T_t^if$ is
continuous $\mu$-almost everywhere on $X$ for every $f\in L^0(X)$ and  $i\in[k]$. Further assume that, for every $i\in[k]$, we have
\begin{align}
\label{oschyp}
\sup_{J \in \NN} \sup_{ I \in \mathfrak S_{J}(\II)}
\norm{ O^r_{I, J} ( T_{t}^{i} f : t \in \II) }_{L^p(X)} 
& \lesssim_{p, r} 
\norm{ f }_{L^p(X)}, 
\qquad f \in L^p(X), 
\end{align}
and
\begin{align}
\label{vvhyp}
\norm[\Big]{ \Big( \sum_{j \in \ZZ} \big( \sup_{t \in \II} \abs{T_t^i} f_j \big)^r \Big)^{1/r} }_{L^{p}(X)}
& \lesssim_{p, r}
\norm[\Big]{ \Big( \sum_{j \in \ZZ} \abs{f_j}^r  \Big)^{1/r} }_{L^{p}(X)}, \qquad (f_j)_{j \in \ZZ} \in L^{p}(X;\ell^r(\ZZ)).
\end{align}
Then we have the following multi-parameter $r$-oscillation estimate:
\begin{align*}
\sup_{J \in \NN} \sup_{ I \in \mathfrak S_{J}(\II^k) }
\norm{ O^r_{I, J} ( T_{t} f : t \in \II^k) }_{L^p(X)} 
\lesssim
\norm{ f }_{L^p(X)}, 
\qquad f \in L^p(X).
\end{align*}
\end{proposition}
In the $\calM_t^{\calP,k',k''} = \calA_t^{\calP,k',k''}$ case, \eqref{vvest} gives us 
\begin{equation}
\label{vvaverage}
 \big\|\sup_{t > 0}|\calA_t^{\calP,k',k''}| f_\frkI\big\|_{L^p(X;\ell^2)} = \big\|\sup_{t > 0}\calA_t^{\calP,k',k''} |f_\frkI|\big\|_{L^p(X;\ell^2)} \lesssim \||f_\frkI|\|_{L^p(X;\ell^2)} = \|f_\frkI\|_{L^p(X;\ell^2)}
\end{equation}
which corresponds to condition \eqref{vvhyp} in the $r=2$ case. We also recall the variation, jump, and one-parameter oscillation inequalities for $\calA_t^{\calP,k',k''}$ and $\calH_t^{\calP,k',k''}$.  
\begin{proposition} \cite[Theorem C]{Troj} \cite[Theorem 1]{MeS} \label{propTMeS}
Let $d, k\ge 1$, $r \in (2,\infty)$, and let $\calP$ be a polynomial mapping as in ~\eqref{polymap}. Let $k',k''\in\{0,1,\ldots,k\}$ with $k'+k''=k$ and let $\calM_t^{\calP,k',k''}$ be either $\calA_t^{\calP,k',k''}$ or $\calH_t^{\calP,k',k''}$. Then, for any $p\in(1,\infty)$, there is a constant $C_{p,d,k,\deg\calP}>0$ such that
\begin{align}
    \norm[\big]{V^r(\calM_t^{\calP,k',k''} f:t>0)}_{L^p(X,\mu)}&\le \frac{r}{r-2}C_{p,d,k,\deg \calP}\|f\|_{L^p(X,\mu)},\label{varest}\\
    \sup_{\lambda>0}\norm[\big]{\lambda N_{\lambda}(\calM_t^{\calP,k',k''} f:t>0)^{1/2}}_{L^p(X,\mu)}&\le C_{p,d,k,\deg \calP}\|f\|_{L^p(X,\mu)},\label{jumpest}\\
    \sup_{N\in\NN}\sup_{I\in\mathfrak{S}_N(\RR_+)}\norm[\big]{O_{I,N}^2(\calM_t^{\calP,k',k''} f:t>0)}_{L^p(X,\mu)}&\le C_{p,d,k,\deg\calP}\norm{f}_{L^p(X,\mu)},\label{oscillest}
\end{align}
for any $f\in L^p(X,\mu)$. The constant $C_{p,d,k,\deg\calP}$ is independent of the
coefficients of the polynomial mapping $\calP$.
\end{proposition}
In particular, \eqref{oscillest} corresponds to condition \eqref{oschyp}. As such, we have the following applications.

\begin{corollary} \label{multiparametercorollary}
    Let $M \in \NN$ and let $(X,\calB,\mu)$ be a  $\sigma$-finite measure space endowed with a family of invertible commuting and measure preserving transformations $S_1^1,\ldots, S_{d_1}^1, \ldots, S_1^M, \ldots, S_{d_M}^M:X\to X$. For each $j \in \{1,\ldots,M\}$, let $\Omega^j$ be a bounded convex open subset of $\RR^{k_j}$ such that $B(0,c) \subseteq \Omega^j \subseteq B(0,1)$ for some $c\in(0, 1)$, let $\calP^j$ be a polynomial mapping
\begin{equation*}
    \calP^j=(\calP_1^j,\dots,\calP_{d_j}^j)\colon\ZZ^{k_j}\to\ZZ^{d_j}
\end{equation*}
where each $\calP_i^j\colon\ZZ^{k_j}\to\ZZ$ is a polynomial of $k_j$ variables with integer coefficients such that $\calP_i^j(0)=0$, and let $k_j',k_j''\in\{0,\ldots,k_j\}$ with $k_j=k_j'+k_j''$. For $f \in L^\infty(X,\mu)$, we define the associated ergodic averages by
\begin{equation*}
\calA_{t}^{\calP^j,k_j',k_j''}f(x):=\frac{1}{\vartheta_{\Omega^j}(t)}\sum_{(n,p)\in\ZZ^{k_j'}\times(\pm\PP)^{k_j''}}f((S_1^j)^{\calP_1^j(n,p)}\cdots (S_{d_j}^j)^{\calP_{d_j}^j(n,p)}x)\mathds{1}_{\Omega_t^j}(n,p)\Big(\prod_{i=1}^{k_j''}\log|p_i|\Big),\quad x\in X,
\end{equation*}
where 
\begin{equation*}
\vartheta_{\Omega^j}(t):=\sum_{(n,p)\in\ZZ^{k_j'}\times(\pm \PP)^{k_j''}}\mathds{1}_{\Omega_t^j}(n,p)\Big(\prod_{i=1}^{k_j''}\log|p_i|\Big).
\end{equation*}
Letting $k = k_1+\ldots+k_M$, $k' = k_1'+\ldots+k_M'$, and $k'' = k_1''+\ldots+k_M''$, we let $(n,p) \in \ZZ^{k'} \times (\pm \PP)^{k''}$ denote $(n_1,p_1,\ldots,n_M,p_M) \in \ZZ^{k_1'} \times (\pm \PP)^{k_1''} \times \ldots \times \ZZ^{k_M'} \times (\pm \PP)^{k_M''} \cong \ZZ^{k'} \times (\pm \PP)^{k''}$. For $f\in L^\infty(X,\mu)$ and $t = (t_1,\ldots,t_M) \in \RR_+^M$, we define the associated multi-parameter ergodic averages by 
\begin{align*}
\calA_{t}f(x) \coloneqq & \calA_{t_1,\ldots,t_M}^{\calP^1, \ldots, \calP^M,k_1',k_1'',\ldots,k_M',k_M''}f(x) \coloneqq  \calA_{t_1}^{\calP^1,k_1',k_1''} \circ \cdots \circ \calA_{t_M}^{\calP^M,k_M',k_M''} f(x)
\\ =\frac{1}{\vartheta(t)}\sum_{(n,p)\in\ZZ^{k'}\times(\pm\PP)^{k''}}  &f((S_1^1)^{\calP_1^1(n_1,p_1)} \cdots (S_{d_1}^1)^{\calP_{d_1}^1(n_1,p_1)} \cdots (S_1^M)^{\calP_1^M(n_M,p_M)}\cdots (S_{d_M}^M)^{\calP_{d_M}^M(n_M,p_M)} x)
\\ &\times \mathds{1}_{\Omega^1_{t_1} \times \ldots \times \Omega^M_{t_M}}(n,p)\Big(\prod_{i=1}^{k''}\log|p_i|\Big),\quad x\in X,
\end{align*}
where 
\begin{equation*}
\vartheta(t):=\sum_{(n,p)\in\ZZ^{k'}\times(\pm \PP)^{k''}}\mathds{1}_{\Omega^1_{t_1} \times \ldots \times \Omega^M_{t_M}}(n,p)\Big(\prod_{i=1}^{k''}\log|p_i|\Big).
\end{equation*}
Let $p\in(1,\infty)$ and $f\in L^p(X,\mu)$. Then we have:
\begin{itemize}
\item[(i)] \textit{(Mean ergodic theorem)} 
the averages $\calA_{t}f$ converge in $L^p(X,\mu)$ norm as $\min\{t_1,\ldots,t_M\} \to \infty$ \label{thm:mean};

\item[(ii)] \textit{(Pointwise ergodic theorem)} 
the averages $\calA_{t}f$ converge pointwise $\mu$-almost everywhere on $X$ as $\min\{t_1,\ldots,t_M\} \to \infty$;

\item[(iii)] \textit{(Maximal ergodic theorem)}
the following maximal estimate holds, including with $p=\infty$:
\begin{align}
\label{thm:max}
\big\|\sup_{t \in \RR_+^M}|\calA_{t}f|\big\|_{L^p(X, \mu)}\lesssim_{d,k,p,M, \deg \mathcal P}\|f\|_{L^p(X, \mu)};
\end{align}
\item[(iv)] \textit{(Oscillation ergodic theorem)}
the following uniform oscillation inequality holds:
\begin{align}
\label{thm:osc}
\sup_{N \in \NN} \sup_{I \in \mathfrak S_{N}(\RR_+^M)}
\norm{O^2_{I,N} (\calA_t f : t \in \RR_+^M)}_{L^p(X)} 
\lesssim_{d,k,p,M, \deg \mathcal P}
\norm{ f }_{L^p(X)}, 
\qquad f \in L^p(X).
\end{align}
\end{itemize}
The implicit constants in \eqref{thm:max} and \eqref{thm:osc} are independent of the
coefficients of the polynomial mapping $\calP$.
\end{corollary}
This generalizes the polynomial Dunford-Zygmund ergodic theorem due to Bourgain as we shall see in section \ref{sec:hist3}. We note that (i) follows from the dominated convergence theorem together with (ii) and (iii), and these each follow from (iv).  Although Corollary \ref{multiparametercorollary} only requires the $\calM_t^{\calP,k',k''} = \calA_t^{\calP,k',k''}$ case of Theorem \ref{MAINthm}, we also prove the $\calM_t^{\calP,k',k''} = \calH_t^{\calP,k',k''}$ case for the sake of independent interest and to exhibit a unified approach that illustrates what common features of the operators are needed in the proof. 

\subsection{Historical background: Classical results and Bourgain.}\label{sec:hist1} 

In 1931, Birkhoff \cite{Birk} and von Neumann \cite{Neu} proved that the averages
\begin{equation}\label{eq:birkhist}
    M_Nf(x):=\frac{1}{N}\sum_{n=1}^Nf(S^nx)
\end{equation}
converge pointwise $\mu$-almost everywhere on $X$ and in $L^p(X,\mu)$ norm respectively for any $f\in L^p(X,\mu)$, $p \in [1,\infty)$, as $N\to\infty$. In 1955, Cotlar \cite{COT} established the pointwise $\mu$-almost everywhere convergence on $X$ as $N\to\infty$ of the ergodic Hilbert transform given by
\begin{equation*}
    H_Nf(x):=\sum_{1\leq|n|\leq N} \frac{f(S^nx)}{n}
\end{equation*}
for any $f\in L^p(X,\mu)$. In 1968, Calder{\'o}n \cite{Cald} made an important observation (now called the Calder{\'o}n transference principle) that some results in ergodic theory can be easily deduced from known results in harmonic analysis. Namely, the convergence of the Birkhoff averages $M_N$ can be deduced from the boundedness of the Hardy--Littlewood maximal function, and the convergence of Cotlar's averages $H_N$ follows from the boundedness of the maximal function for the truncated discrete Hilbert transform. As we will see ahead, this observation has had a huge impact in the study of convergence problems in ergodic theory. 

We briefly sketch the classical approach of handling the problem of pointwise convergence. It consists of two steps:
\begin{itemize}
    \item[(a)] Establish $L^p$-boundedness for the corresponding maximal function.
    \item[(b)] Find a dense class of functions in $L^p(X,\mu)$ for which the pointwise convergence holds. 
\end{itemize}
In the case of Birkhoff's averages $M_N$, the Calder{\'o}n transference principle allows one to deduce the estimate
\begin{equation*}
\|\sup_{N\in\NN}|M_Nf|\|_{L^p(X,\mu)}\lesssim_{p}\|f\|_{L^p(X,\mu)}
\end{equation*}
for $p\in(1,\infty]$ from the estimate for the discrete Hardy--Littlewood maximal function (and we have a weak-type estimate for $p=1$). In turn, estimates for the discrete Hardy--Littlewood maximal function follow easily from those for the continuous one. This establishes the first step (a). For the second step, one can use the idea of Riesz decomposition \cite{Riesz} to analyze the space ${\II}_S\oplus {\TT}_S\subseteq L^2(X,\mu)$, where
\begin{align*}
\qquad\qquad { \II}_S:=\{f\in L^2(X,\mu): f\circ S =f\}
\qquad \text{ and } \qquad
{\TT}_S:=\{h\circ S-h: h\in L^2(X,\mu)\cap L^{\infty}(X,\mu)\}.
\end{align*}
We see that $M_Nf=f$ for  $f\in {\II}_S$ and, for $g=h\circ S - h\in {\TT}_S$, we have 
\begin{equation*}
M_Ng(x)=\frac{1}{N}\big(h(S^{N+1}x)-h(Sx)\big)
\end{equation*}
by telescoping. Consequently, we see that $M_Ng\to0$ as $N\to\infty$.  This establishes $\mu$-almost everywhere pointwise convergence of $M_N$ on ${\II}_S\oplus {\TT}_S$, which is dense
in $L^2(X,\mu)$. Since $L^2(X,\mu)$ is dense in $L^p(X,\mu)$ for every $p\in[1,\infty)$, this establishes (b).

At the beginning of the 1980's, Bellow \cite{Bel} and independently Furstenberg \cite{Fus} posed the problem of the pointwise convergence of the averages along squares given by
\begin{equation*}
    T_N f(x):=\frac{1}{N}\sum_{n=1}^N{f(S^{n^2}x)}.
\end{equation*}
Despite its similarity to Birkhoff's theorem, the problem of pointwise convergence of the $T_N$ averages has a totally different nature from that of its linear counterpart, and the standard approach is insufficient in this case. For the first step, by the Calder{\'o}n transference principle, it is enough to establish $\ell^p$ bounds for the maximal function given by
\begin{equation}\label{eq:maximalsq}
    \sup_{N\in\NN}\frac{1}{N}\sum_{n=1}^N f(x-n^2),\quad f\in\ell^p(\ZZ).
\end{equation}
The $\ell^p$ estimate for the above maximal function does not follow directly from the continuous counterpart and requires completely new methods. However, a more serious problem arises in connection with the second step. The telescoping idea fails in the case of the averages $T_N g$ since the $(n+1)^2-n^2 = 2n+1$ gap sizes are unbounded. 

At the end of the 1980's, Bourgain established the pointwise convergence of the averages $T_N$ in a series of groundbreaking articles \cite{B1,B2,B3}. By using the Hardy--Littlewood circle method from analytic number theory, he established $\ell^p$-bounds for the maximal function \eqref{eq:maximalsq}, which establishes step (a). He then bypassed the problem of finding the requisite dense class of functions by using the oscillation seminorm. Bourgain \cite{B3} proved that,
for any $\lambda>1$ and any sequence of integers $I=(I_j:{j\in\NN})$ with $I_{j+1}>2I_j$ for all $j\in\NN$, we have
\begin{align}\label{eq:non-uni}
\norm[\big]{O_{I, N}^2(T_{\lambda^n}f:n\in\NN)}_{L^2(X,\mu)}
\le C_{I,\lambda}(N)\norm{f}_{L^2(X,\mu)}, \qquad N\in\NN,
\end{align}
for any $f\in L^2(X,\mu)$ with $\lim_{N\to\infty} N^{-1/2}C_{I, \lambda}(N)=0$. Inequality \eqref{eq:non-uni} suffices to establish the pointwise convergence of the averaging operators $T_Nf$ for any  $f\in L^2(X, \mu)$ (see \cite[Proposition 2.8]{MSWr} for details even in the multi-parameter setting). Indeed, it can be thought of as the weakest possible quantitative form of pointwise convergence since one can derive \eqref{eq:non-uni} with $C_{I,\lambda}(N)$ at most $N^{1/2}$ from the $\ell^2$ bound for the maximal function \eqref{eq:maximalsq}.

In the same series of papers, by similar methods, Bourgain established the pointwise convergence of the averages along primes
\begin{equation*}
    \frac{1}{|\PP_N|}\sum_{n=1}^Nf(S^nx)\mathds{1}_{\PP}(n)
\end{equation*}
for $f\in L^p(X,\mu)$ with $p>\frac12(1+\sqrt{3})$. In the same year, Wierdl \cite{Wie} extended Bourgain's result to $p\in(1,\infty)$.

\subsection{Historical background: Developments towards Proposition \ref{propTMeS}}\label{sec:hist2}

An intriguing question was the issue of uniformity in the inequality \eqref{eq:non-uni}. Shortly after the groundbreaking work of Bourgain, Lacey \cite[Theorem 4.23, p. 95]{RW} improved inequality \eqref{eq:non-uni} showing that, for every $\lambda>1$,
there is a constant $C_{\lambda}>0$ such that
\begin{align}
\label{eq:Lacey}
\sup_{N\in\NN}\sup_{I\in \mathfrak S_N(\LL_{\tau})}\norm[\big]{O_{I, N}^2(T_{\lambda^n}f:n\in\NN)}_{L^2(X)}
\le C_{\lambda}\norm{f}_{L^2(X)},\quad f\in L^2(X,\mu),
\end{align}
where $\LL_{\tau}:=\{\tau^n:n\in \NN\}$. This result motivated the question about uniform estimates independent of $\lambda>1$ in \eqref{eq:Lacey}. In the case of Birkhoff's averages, this is the Rosenblatt-Wierdl conjecture explicitly formulated in \cite[Problem 4.12, p. 80]{RW} in the early 1990s.

In 1998, Jones, Kaufman, Rosenblatt, and
Wierdl \cite{jkrw} established the uniform oscillation inequality on $L^p(X,\mu)$ for the standard Birkhoff averages $M_N$. Two years later, Campbell, Jones, Reinhold, and Wierdl \cite{CJRW1} established the uniform  oscillation inequality for the ergodic Hilbert transform. In 2003, Jones, Rosenblatt, and Wierdl \cite{jrw} proved uniform oscillation inequalities on $L^p(X,\mu)$ with $p\in(1,2]$ for the Birkhoff averages over cubes. However, the case of polynomial averages, even one-dimensional, was open until recent works \cite{MSS,S}, see \cite{MeS} for the case of averages over primes.

In 2015, Mirek and Trojan \cite{MT}, using the ideas of Bourgain and Wierdl, established $\mu$-almost everywhere pointwise convergence of the Cotlar averages along the primes, 
\begin{equation*}
    \sum_{p\in(\pm\PP_N)} \frac{f(S^p)}{p}\log |p|.
\end{equation*}
They proved that the corresponding maximal function is bounded on $L^p(X,\mu)$ with $p>1$ and showed that the analogue of Bourgain's non-uniform oscillation inequality \eqref{eq:non-uni} holds for those averages.

In the same year, Zorin-Kranich \cite{zk} established the pointwise convergence of the averages related to the polynomial mapping given by
\begin{equation*}
    \Tilde{\calP}=(n,n^2,n^3,\ldots,n^d)\colon\ZZ\to\ZZ^d.
\end{equation*}
Namely, he proved that, for any $r>2$ and $|\frac{1}{p}-\frac{1}{2}|<\frac{1}{2(d+1)}$, we have the following $r$-variational estimate
\begin{equation*}
\|V^r(\calA_N^{\Tilde{\calP},1,0}f:N\in\NN)\|_{L^p(X,\mu)}\lesssim_{p,r}\|f\|_{L^p(X,\mu)}.
\end{equation*}
As a consequence, the averages $\calA_N^{\Tilde{\calP},1,0}f$ converge $\mu$-almost everywhere for any $f\in L^p(X,\mu)$. 

In 2016, Mirek and Trojan \cite{MT1} established the pointwise convergence for the averages \eqref{def:average} taken over cubes with $k'=k$, that is
\begin{equation*}
    A_{N,{\rm cube}}^{\calP,k,0} f(x):=\frac{1}{N^k}\sum_{y\in[0,N]^k\cap\ZZ^k}f(S_1^{\calP_1(y)}S_2^{\calP_1(y)}\cdots S_d^{\calP_d(y)}x).
\end{equation*}
There, Mirek and Trojan noted for the first time that the Rademacher--Menshov inequality~\eqref{eq:remark3} may be used to establish $r$-variational estimates. For $p\in(1,\infty)$ and $r>\max\{p, p/(p-1)\}$, they proved that
\begin{equation*}
    \|V^r(A_{N,{\rm cube}}^{\calP,k,0} f:N\in\NN)\|_{L^p(X,\mu)}\leq  C_{p,d,k,{\rm deg}\calP}\|f\|_{L^p(X,\mu)}.
\end{equation*}

Unfortunately, the methods introduced by Bourgain had limitations. These work perfectly fine in the case of the $L^2$ estimates, but, in the case of an $L^p$ estimates with $p\neq2$, there arise difficulties which are hard to overcome concerning the fractions around which major arcs are defined. However, Ionescu and Wainger \cite{IW}, in their groundbreaking 2005 work about discrete singular Radon operators, introduced a set of fractions for which the circle method can be applied towards $L^p$ estimates with $p\neq2$. 

In 2015, Mirek \cite{M} built a discrete counterpart of the Littlewood--Paley theory using the Ionescu--Wainger multiplier theorem and used it to reprove the main result from \cite{IW}. In 2017, Mirek, Stein, and Trojan \cite{MST1,MST2} further exploited these ideas together with the Rademacher--Menshov inequality from \cite{MT1} to obtain an $L^p$ estimate for the $r$-variation seminorm for both $\calA_t^{\calP,k,0}$ and $\calH_t^{\calP,k,0}$ associated with convex sets in the full range of parameters. Namely, they showed that
\begin{equation}\label{eq:var:MST}
\big\|V^r(\calM_t^{\calP,k,0}f: t>0)\big\|_{L^p(X, \mu)}\lesssim_{d,k,p, r,  \deg \mathcal P}\|f\|_{L^p(X, \mu)}
\end{equation}
for $p\in(1,\infty)$ and $r\in(2,\infty)$, where $\calM_t^{\calP,k,0}$ is either $\calA_t^{\calP,k,0}$ or $\calH_t^{\calP,k,0}$.

In 2019, Trojan \cite{Troj} proved the variation case of Proposition \ref{propTMeS}. A straightforward consequence is the $\mu$-almost everywhere convergence of the averages $\calA_t^{\calP,k',k''}f$ and $\calH_t^{\calP,k',k''}f$. 

In 2020, Mirek, Stein, and Zorin-Kranich \cite{MSZ3} further refined the methods developed in \cite{MST1,MST2} and proved a uniform $L^p$ estimate for the $\lambda$-jump counting function. They proved that  
\begin{equation}\label{eq:jump:MSZK}
    \sup_{\lambda>0}\norm[\big]{\lambda N_{\lambda}(\calM_t^{\calP,k,0}f:t>0)^{1/2}}_{L^p(X,\mu)}\le C_{p,d,k,{\rm deg}\calP}\|f\|_{L^p(X,\mu)}
\end{equation}
for any $p\in(1,\infty)$ and any $f\in L^p(X,\mu)$, where $\calM_t^{\calP,k,0}f$ is either $\calA_t^{\calP,k,0}f$ or $\calH_t^{\calP,k,0}f$. There, the operators $\calH_t^{\calP,k,0}$ are associated with Calder{\'o}n--Zygmund kernels satisfying the H{\"o}lder continuity condition generalizing \eqref{eq:K-modulus-cont}:
For some $\sigma \in (0,1]$ and for every $x, y\in\RR^k\setminus\{0\}$ with $2|y|\leq |x|$, we have
\begin{equation}
\label{eq:K-Holder}
\abs{K(x)-K(x+y)}
\lesssim 
\abs{y}^{\sigma} \abs{x}^{-(k+\sigma)}.
\end{equation}
It is worth noting that the inequality \eqref{eq:jump:MSZK} implies the $r$-variation inequality \eqref{eq:var:MST}. 

In 2021, Mirek, Słomian, and Szarek~\cite{MSS} established the oscillation inequality
\begin{equation}\label{osc:MSS}
\sup_{N\in\NN}\sup_{I\in\mathfrak{S}_N(\RR_+)}\norm[\big]{O_{I,N}^2(\calA_t^{\calP,k,0}f:t>0)}_{L^p(X,\mu)}\le C_{p,d,k,{\rm deg}\calP}\norm{f}_{L^p(X,\mu)},
\end{equation}
and, Słomian \cite{S} later proved the counterpart of \eqref{osc:MSS} in the case of the operators $\calH_t^{\calP,k,0}$ related to Calder{\'o}n--Zygmund kernels satisfying \eqref{eq:K-Holder}. Lastly, in 2023, Słomian and the author \cite{MeS} proved the oscillation and jump cases of Proposition \ref{propTMeS}. This is the most general quantitative version of the one parameter ergodic theorem for both averages $\calA_t^{\calP,k',k''}$ and $\calH_t^{\calP,k',k''}$ (cf. \cite[Theorem 1.20]{MISZWR}), which concludes the work of many authors over the last decades.

\subsection{Historical background: Multi-parameter problems.}\label{sec:hist3} 

In 1951, Dunford \cite{D} and independently
Zygmund \cite{Z} showed that the two-step procedure can be applied in a multi-parameter setting. Even for $S_1,\ldots,S_d$ not necessarily commuting, the Dunford--Zygmund ergodic theorem states that the averages 
\begin{align*}
A_{M_1,\ldots, M_d;S_1,\ldots,S_d}^{n_1,\ldots, n_d}f(x)
\coloneqq \frac{1}{M_1\cdots M_d} \sum_{n_1=1}^{M_1} \ldots \sum_{n_d=1}^{M_d} f(S_1^{n_1}\cdots S_d^{n_d}x), \qquad x \in X, 
\end{align*}
converge almost everywhere on $X$ and in $L^p(X)$ norm as
$\min\{M_1,\ldots, M_d\}\to\infty$ for every $f\in L^p(X)$, $p\in(1, \infty)$. Using the identity
\begin{align*}
A_{M_1,\ldots, M_d;S_1,\ldots,S_d}^{n_1,\ldots, n_d}f =
A_{M_1;S_1}^{n_1} \circ \cdots \circ A_{M_d;S_d}^{n_d}f,
\end{align*}
the $L^p(X)$, $p\in(1, \infty]$, bounds for the strong
maximal function
$\sup_{M\in\NN^d}|A_{M_1,\ldots, M_d;S_1,\ldots,S_d}^{n_1,\ldots, n_d}f|$
follow by applying $d$ times the corresponding $L^p(X)$ bounds
for $\sup_{M\in\NN}|A_{M;S}^{n}f|$. This establishes (a), and (b) can be established by a suitable adaptation of the telescoping argument to the
multi-parameter setting and an application of the classical Birkhoff
ergodic theorem, see \cite{Nevo} for more details. We note that that the operator
$f \mapsto \sup_{M\in \NN^d}|A_{M_1,\ldots, M_d;S_1,\ldots,S_d}^{n_1,\ldots, n_d}f|$
is not of weak type $(1, 1)$ in general, so the pointwise convergence may fail if $p = 1$. A model example is $X=\mathbb{Z}^d$ with $S_j x = x - e_j$, $1 \leq j \leq d$, where $e_j$ is the $j$th coordinate vector. It is well known that the weak type $(1,1)$ estimate does not hold for the corresponding strong maximal operator, see \cite[Section X.2.3]{bigs}. 

After completing \cite{B1,B2,B3}, Bourgain observed that the
Dunford--Zygmund ergodic theorem can be extended to the
polynomial setting at the expense of imposing that the 
measure-preserving transformations commute. Bourgain's result can be formulated as follows. 

\begin{proposition}[Polynomial Dunford--Zygmund ergodic theorem]\cite[Theorem 1.25]{MISZWR} \label{PDZ}
    Let $M \in \NN$, let $(X,\calB,\mu)$ be a  $\sigma$-finite measure space endowed with a family of invertible commuting and measure preserving transformations $S_1,\ldots, S_M:X\to X$, and consider a polynomial mapping
\begin{equation*}
    \calP=(\calP_1,\dots,\calP_{M})\colon \ZZ^M \to \ZZ^M
\end{equation*}
where each $\calP_j\colon\ZZ \to \ZZ$ is a polynomial of one variable with integer coefficients such that $\calP_j(0)=0$. For $f \in L^\infty(X,\mu)$ and $t,t_1,\ldots,t_M \in \NN$, we define the associated ergodic averages by
\begin{equation*}
A_{t}^{\calP_j}f(x) \coloneqq \frac{1}{t}\sum_{n=1}^t f(S_1^{\calP_j(n)} x), \quad x\in X,
\end{equation*}
and  
\begin{align*}
A_{t_1, \ldots, t_M}^{\calP_1, \ldots, \calP_M} f(x) \coloneqq  A_{t_1}^{\calP^1} \circ \cdots \circ \calA_{t_M}^{\calP_M,} f(x) =\frac{1}{t_1\cdots t_M}\sum_{n_1=1}^{t_1} \ldots \sum_{n_M=1}^{t_M} f(S_1^{\calP_1(n_1)} \cdots S_{M}^{\calP_M(n_M)} x), \quad x\in X.
\end{align*}
Let $p\in(1,\infty)$ and $f\in L^p(X,\mu)$. Then we have:
\begin{itemize}
\item[(i)] \textit{(Mean ergodic theorem)} 
the averages $A_{t_1, \ldots, t_M}^{\calP_1, \ldots, \calP_M} f$ converge in $L^p(X,\mu)$ norm as $\min\{t_1,\ldots,t_M\} \to \infty$ ;

\item[(ii)] \textit{(Pointwise ergodic theorem)} 
the averages $A_{t_1, \ldots, t_M}^{\calP_1, \ldots, \calP_M} f$ converge pointwise $\mu$-almost everywhere on $X$ as $\min\{t_1,\ldots,t_M\} \to \infty$;

\item[(iii)] \textit{(Maximal ergodic theorem)}
the following maximal estimate holds, including with $p=\infty$:
\begin{align*}
\big\|\sup_{t \in \NN^M}|A_{t_1, \ldots, t_M}^{\calP_1, \ldots, \calP_M} f|\big\|_{L^p(X, \mu)}\lesssim_{p,M, \deg \mathcal P}\|f\|_{L^p(X, \mu)};
\end{align*}
\item[(iv)] \textit{(Oscillation ergodic theorem)}
the following uniform oscillation inequality holds:
\begin{align*}
\sup_{N \in \NN} \sup_{I \in \mathfrak S_{N}(\RR_+^M)}
\norm{O^2_{I,N} (A_{t_1, \ldots, t_M}^{\calP_1, \ldots, \calP_M} f : t \in \NN^M)}_{L^p(X)} 
\lesssim_{p,M, \deg \mathcal P}
\norm{ f }_{L^p(X)}, 
\qquad f \in L^p(X),
\end{align*}
\end{itemize}
with implicit constants independent of the
coefficients of the polynomial mapping $\calP$.
\end{proposition}

Proposition \ref{PDZ}(i)-(iii) is attributed to Bourgain, though it was never published, see \cite{MISZWR} for a proof and additional historical notes, and Proposition \ref{PDZ}(iv) with linear polynomials $\calP_1(t) = \ldots = \calP_M(t) = t$ was established in \cite{jrwRect}. Corollary \ref{multiparametercorollary} is a significant generalization of Proposition \ref{PDZ} (one may check that the proof is easily adaptable to sums taken over $\NN$ instead of $\ZZ$) in that it allows for averages taken over primes and over more general polynomial orbits. Indeed, from the point of view of the permitted polynomial orbits, Corollary \ref{multiparametercorollary} is the most one can extend Proposition \ref{PDZ} without having to go beyond averaging operators that can be written as the composition of one-parameter averaging operators. For comparison, proving the analogue of Proposition \ref{PDZ} for averages of the form 
\begin{align*}
\frac{1}{t_1\cdots t_M}\sum_{n_1=1}^{t_1} \ldots \sum_{n_M=1}^{t_M} f(S_1^{\calP_1(n_1,\ldots,n_M)} \cdots S_{M}^{\calP_M(n_1,\ldots,n_M)} x)
\end{align*}
is a central open problem in modern ergodic theory that can be seen as a multi-parameter variant of the Bellow and Furstenberg problem (cf. \cite[Conjecture 1.29]{MISZWR}, and see \cite{BMSW} for some recent progress).

\section{Notation and necessary tools}
\subsection{Basic notation}\label{sec:notation}
We denote $\NN:=\{1, 2, \ldots\}$, $\NN_0:=\{0,1,2,\ldots\}$, and $\RR_+:=(0, \infty)$. For
$d\in\NN$, the sets $\ZZ^d$, $\RR^d$, $\CC^d$, and $\TT^d = (\RR/\ZZ)^d \equiv[-1/2, 1/2)^d$ have the standard meanings. For each $N\in\NN$, we set
    \[
    \NN_N:=\{1,\ldots, N\}.
    \]  
For any $x\in\RR$, we set
\[
\lfloor x \rfloor: = \max\{ n \in \ZZ : n \le x \}.
\]
For $u\in\NN$, we define the set
\begin{equation*}
     2^{u\NN}:=\{2^{un}\colon n\in\NN\}.
\end{equation*}

For two non-negative numbers $A$ and $B$, we write $A \lesssim B$ to indicate that $A\le CB$ for some $C>0$ that may change from line to line, and we may write $\lesssim_{\delta}$ if the implicit constant depends on $\delta$.

We denote the standard inner product on $\RR^d$ by $ x\cdot\xi$. Moreover, for any $x\in\RR^d$, we denote the $\ell^2$-norm and the maximum norm respectively by
\begin{align*}
\abs{x}:=\abs{x}_2:=\sqrt{\ipr{x}{x}} \qquad \text{ and } \qquad |x|_{\infty}:=\max_{1\leq k\leq d}|x_k|.
\end{align*}

For a multi-index $\gamma=(\gamma_1,\dots,\gamma_k)\in\N^k_0$, we abuse the notation to write
$|\gamma|:=\gamma_1+\cdots+\gamma_k$. No confusion should arise since all multi-indices will be denoted by $\gamma$.

 Let $\II \subset \RR$ and let $E$ be either of $\RR^d$ or $\ZZ^d$ with the usual measures and let $(f_{\frkI,t}:\frkI \in \NN) \in L^p(E,\ell^2)$ for all $t \in \II$. We define 
\[\calS_E^p(f_{\frkI,t} : t \in \II) \coloneqq \big\|\sup_{t\in \II} |f_{\frkI,t}-f_{\frkI,\inf\II}|\big\|_{L^p(E;\ell^2)} = \bigg\|\bigg(\sum_{\frkI=1}^\infty \sup_{t\in \II} |f_{\frkI,t}-f_{\frkI,\inf\II}|^2\bigg)^{1/2}\bigg\|_{L^p(E)}.\]
As we shall see, working with this rather than the usual maximal function is just a technical adaptation to be more similar to the variation, jump, and oscillation quantities that have been studied before. 

\subsection{Rademacher--Menshov inequality}
We recall a basic numerical inequality. A variational version of this inequality was proven by Lewko–Lewko \cite[Lemma
13]{LL}, see also \cite[Lemma 2.5, p. 534]{MSZ2}. 
\begin{proposition}\label{rem:RM} 
For any $k,m\in\NN$ with $k<2^m$ and any sequence of complex numbers $(a_n:n\in\NN)$, we have
\begin{equation}\label{eq:remark3}
    \sup_{k \leq n\leq 2^m} |a_n| \leq |a_k| + \sqrt{2} \sum_{i=1}^s\Big(\sum_{j}{|a_{u_{j+1}^i}-a_{u_j^i}|^2}\Big)^{1/2},
\end{equation}
where each $[u_j^i,u_{j+1}^i)$ is a dyadic interval contained in $[k, 2^m]$ of the form $[j2^i,(j+1)2^{i})$ for some $0\le i\le m$ and $0\le j\le 2^{m-i}-1$. 
\end{proposition}


\subsection{Marcinkiewicz--Zygmund inequality}
We recall a result extending the Marcinkiewicz–-Zygmund inequality to the Hilbert space setting. Let $(T_m: m \in \NN_0)$ be a family of bounded linear operators, $T_m: L^p(X) \rightarrow L^p(X)$. For each $\omega \in [0,1]$, we define
\begin{equation*}
T^\omega = \sum_{m \in \NN_0} \epsilon_m(\omega)T_m
\end{equation*}
where $(\epsilon_m: m \in \NN_0)$ is the sequence of Rademacher functions on $[0,1]$. 
\begin{proposition}\cite[Lemma 2.1]{MST1} \label{MZprop} Let $p \in (0,\infty)$. Suppose there is a constant $C_p > 0$ such that, for all $\omega \in [0,1]$ and $f \in L^p(X)$, we have
\begin{equation*}
    \|T^\omega f\|_{L^p(X)} \leq C_p \|f\|_{L^p(X)},
\end{equation*}
then there is a constant $C$ such that
\begin{equation}\label{MZgen}
    \big\|\big(\sum_{m \in \NN_0} |T_m f_\frkI|^2\big)^{1/2}\big\|_{L^p(X;\ell^2)} \leq C C_p \|f_\frkI\|_{L^p(X;\ell^2)}
\end{equation}
for every sequence of functions $(f_\frkI)_{\frkI \in \NN}$ in $L^p(X;\ell^2)$. Moreover, if $T_m \equiv 0$ for all $m \in \NN$, then \eqref{MZgen} recovers the Marcinkiewicz--Zygmund inequality 
\begin{equation}\label{MZbasic}
    \|T_0 f_\frkI\|_{L^p(X;\ell^2)} \leq C C_p \|f_\frkI\|_{L^p(X;\ell^2)}.
\end{equation}
\end{proposition}

\subsection{Reductions: Calder{\'o}n transference and lifting}\label{sec:red}
By the Calder\'on transference principle \cite{Cald}, we may restrict attention to the model dynamical system of $\ZZ^d$ equipped with the counting measure and the shift operators $S_j\colon\ZZ^d\to\ZZ^d$ given by $S_j(x_1,\ldots,x_d):=(x_1,\ldots,x_j-1,\ldots,x_d)$. We denote the corresponding averaging operators by
\[ 
A_t^{\calP,k',k''} f(x) = 
	\frac{1}{\vartheta_\Omega(t)}
	\sum_{(n,p)\in\ZZ^{k'}\times(\pm\PP)^{k''}}
	f\big(x - \calP(n, p)\big)\ind{\Omega_t}(n, p) 
	\Big(\prod_{j = 1}^{k''} \log |p_j|\Big) 
\]
and
\[ 
H_t^{\calP,k',k''} f(x) = \sum_{(n,p)\in\ZZ^{k'}\times(\pm\PP)^{k''}}
	f\big(x - \calP(n, p)\big) K(n, p) \ind{\Omega_t}(n, p) \Big(\prod_{j = 1}^{k''} \log \abs{p_j}\Big). 
\]
Moreover, by a standard lifting argument, it suffices to
prove Theorem~\ref{MAINthm} for a canonical case of the polynomial mapping $\calP$. Let $\mathcal{P}$ be a polynomial mapping as in \eqref{polymap}. We define 
\begin{equation*}
    {\rm deg}\, \mathcal{P}:=\max\{{\rm deg}\, \mathcal{P}_j: 1\le j\le d\}
\end{equation*}
and consider the set of multi-indices 
\begin{equation*}
    \Gamma:=\big\{\gamma\in\N_0^k\setminus\{0\}: 0<|\gamma|\le {\rm deg}\, \mathcal{P}\big\}
\end{equation*}
equipped with the lexicographic order. 
We define the \emph{canonical polynomial mapping} by
\begin{equation}\label{canonicalpolymap}
    \RR^k\ni x=(x_1,\dots,x_k)\mapsto\calQ(x):=(x^\gamma\colon\gamma\in\Gamma)\in\RR^\Gamma,
\end{equation} 
where $x^\gamma=x_1^{\gamma_1}x_2^{\gamma_2}\cdots x_k^{\gamma_k}$.
By invoking the lifting procedure described in \cite[Lemma 2.2]{MST1} (see also \cite[Chapter XI]{bigs}), the following implies Theorem \ref{MAINthm}.
\begin{theorem}\label{maintheorem2}
Let $k\in\NN$, let $\Gamma \subset \NN^{k} \setminus \{0\}$ be a nonempty finite set, and let $k',k''\in\{0,1,\ldots,k\}$ with $k'+k''=k$. Let $M_t^{k',k''}$ be either $A_t^{\calQ,k',k''}$ or $H_t^{\calQ,k',k''}$. For any $p\in(1,\infty)$, there is a constant $C_{p,k,|\Gamma|}>0$ such that 
\begin{equation}
        \big\|\sup_{t > 0}|M_t^{\calP,k',k''} f_\frkI|\big\|_{\ell^p(\ZZ^{\Gamma};\ell^2)} \leq C_{p,k,|\Gamma|}\|f_\frkI\|_{\ell^p(\ZZ^{\Gamma};\ell^2)}.
\end{equation}
\end{theorem}
\subsection{Fourier transform and Ionescu--Wainger multiplier theorem} 

Let $\GG=\RR^d$ or $\GG=\ZZ^d$ and let $\GG^*$ denote the dual group of $\GG$. For every $z\in\CC$, we set $\ex(z):=e^{2\pi {\bm i} z}$, where ${\bm i}^2=-1$. Let $\calF_{\GG}$ denote the Fourier transform on $\GG$ defined for any $f \in L^1(\GG)$ by
\begin{align*}
\calF_{\GG} f(\xi) := \int_{\GG} f(x) \ex(x\cdot\xi) {\rm d}\mu(x),\quad \xi\in\GG^*,
\end{align*}
where $\mu$ is the usual Haar measure on $\GG$. For any bounded function $\mathfrak m\colon \GG^*\to\CC$, we define the corresponding Fourier multiplier operator by 
\begin{align}
\label{eq:fourier mult}
T_{\GG}[\mathfrak m]f(x):=\int_{\GG^*}\ex(-\xi\cdot x)\mathfrak m(\xi)\calF_{\GG}f(\xi){\rm d}\xi, \quad x\in\GG.
\end{align}
Here, we assume that $f\colon\GG\to\CC$ is a compactly supported function on $\GG$ (and smooth if $\GG=\RR^d$) or any other function for which \eqref{eq:fourier mult} makes sense.

An indispensable tool in the proof of Theorem~\ref{maintheorem2} is the vector-valued Ionescu--Wainger multiplier theorem from \cite[Section 2]{MSZ3} with an improvement by Tao \cite{TaoIW}.
\begin{theorem}\label{IW}
For every $\varrho>0$, there exists a family $(P_{\leq N})_{N\in\NN}$ of subsets of $\NN$ such that:  
\begin{enumerate}[label*={(\roman*)}]
\item \label{IW1} $\NN_N\subseteq P_{\leq N}\subseteq\NN_{\max\{N, e^{N^{\varrho}}\}}$.
\item \label{IW2}  If $N_1\le N_2$, then $P_{\leq N_1}\subseteq P_{\leq N_2}$.
\item \label{IW3}  If $q \in P_{\leq N}$, then all factors of $q$ also lie in $P_{\leq N}$.
\item \label{IW4} $\lcm(P_N) \leq 3^N$. 
\end{enumerate}

Furthermore, for every $p \in (1,\infty)$, there exists  $0<C_{p, \varrho, |\Gamma|}<\infty$ such that, for every $N\in\NN$, the following holds:

Let $0<\varepsilon_N \le e^{-N^{2\varrho}}$ and let $\mathbf Q:=[-1/2, 1/2)^\Gamma$ be a unit cube. Let $\frkM \colon\RR^{\Gamma} \to L(H_0,H_1)$ be a measurable function supported on $\varepsilon_{N}\mathbf Q$ taking values in $L(H_{0},H_{1})$, the space  of bounded linear operators between separable Hilbert spaces $H_{0}$ and $H_{1}$.
Let $0 \leq \mathbf A_{p} \leq \infty$ denote the smallest constant such that
\[
\big\|T_{\RR^\Gamma}[\frkM ]f\big\|_{L^{p}(\RR^{\Gamma};H_1)}
\leq
\mathbf A_{p} \|f\|_{L^{p}(\RR^{\Gamma};H_0)}
\]
for every function $f\in L^2(\RR^\Gamma;H_0)\cap L^{p}(\RR^\Gamma;H_0)$. Then, the multiplier
\[
\Delta_N(\xi)
:=\sum_{b \in\Sigma_{\leq N}}
\frkM (\xi - b),
\]
where $\Sigma_{\leq N}$ is defined by
\[
\Sigma_{\leq N} := \Big\{ \frac{a}{q}\in\QQ^\Gamma\cap\TT^\Gamma:  q \in P_{\leq N}\text{ and } {\rm gcd}(a, q)=1\Big\},
\]
satisfies
\begin{equation}\label{eq:iw:main}
\big\|T_{\ZZ^\Gamma}[\Delta_{N}]f\big\|_{\ell^p(\ZZ^{\Gamma};H_1)}
\le C_{p,\varrho,|\Gamma|} (\log N) \mathbf A_{p}
\|f\|_{\ell^p(\ZZ^{\Gamma};H_0)}
\end{equation}
for every $f\in \ell^p(\ZZ^\Gamma;H_0)$, (cf. \cite[Theorem 1.4]{TaoIW} which removes the factor of $\log N$ in the inequality~\eqref{eq:iw:main}).
\end{theorem}

\subsection{Exponential sums} In this section, we present some general results concerning the behavior of exponential sums. The following proposition is an enhancement of the variant of Weyl's inequality due to Trojan \cite[Theorem 2]{Troj} that allows us to estimate exponential sums related to a possibly non-differentiable function $\phi$, (cf. \cite[Theorem A.1]{MSZ3}).
\begin{proposition}[Weyl's inequality] \cite[Proposition 6]{MeS}\label{weylinq}
Let $\alpha >0$, $k\in\NN$, and let $\Gamma \subset \NN^{k} \setminus \{0\}$ be a nonempty finite set. Let $\Omega'\subseteq\Omega\subseteq B(0,N)\subset\RR^k$ be convex sets and let $\phi\colon \Omega\cap\ZZ^{k}\to\CC$. There is $\beta_\alpha>0$ such that, for any $\beta > \beta_\alpha$, if there is a multi-index $\gamma_0\in\Gamma$ with
\begin{align*}
\abs[\Big]{\xi_{\gamma_0} - \frac{a}{q}}
\leq
\frac{1}{q^2}
\end{align*}
for some coprime integers $a$ and $q$ with $1\leq a\leq q$ and $(\log N)^\beta\leq q\leq N^{|\gamma_0|}(\log N)^{-\beta}$, then
\begin{align*}
\abs[\Big]{\sum_{(n,p)\in\ZZ^{k'}\times(\pm\PP)^{k''}}\ex(\xi\cdot\calQ(n,p))\phi(n,p)\mathds{1}_{\Omega\setminus\Omega'}(n,p)}
&\lesssim
N^k\log(N)^{-\alpha}\norm{\phi}_{L^\infty(\Omega\setminus\Omega')} \\
&\,\,+N^{k} \sup_{\substack{\abs{x-y}\leq N(\log N)^{-\alpha}\\x,y\in\Omega\setminus\Omega'}} \abs{\phi(x)-\phi(y)}.
\end{align*}
The implicit constant is independent of the function $\phi$, the variable $\xi$, the sets $\Omega,\Omega'$, and the numbers $a$, $q$, and $N$.
\end{proposition}

The next result is a generalization of \cite[Proposition 4.1]{Troj} and \cite[Proposition 4.2]{Troj} in the spirit of \cite[Proposition 4.18]{MSZ3}. For $q\in\NN$ and $a\in\NN_q^\Gamma$ with ${\rm gcd}(a,q)=1$, the \emph{Gaussian sum} related to the polynomial mapping $\calQ$ is given by
\begin{equation}
G(a/q) := \frac{1}{q^{k'}} \frac{1}{\varphi(q)^{k''}} \sum_{x \in \NN^{k'}_q} \sum_{y \in A_q^{k''}} 
	\ex((a/q)\cdot\calQ(x, y)),
\end{equation}
where $A_q:=\{a\in\NN_q:{\rm gdc}(a,q)=1\}$ and $\varphi$ is Euler's totient function. There is $\delta > 0$ such that
\begin{equation}\label{gausssumest}
    \big|G(a/q) \big| \lesssim q^{-\delta},
\end{equation}
according to \cite[Theorem 3]{Troj}.

\begin{proposition} \cite[Lemma 7]{MeS}\label{approxlemma}
Let $N\in\NN$ and let $\Omega\subseteq B(0,N)\subset\RR^k$ be a convex set or a Boolean combination
of finitely many convex sets. Let $\calK\colon\RR^k\to\CC$ be a continuous function supported in $\Omega$. Then, for each $\beta>0$, there is a constant $c = c_{\beta}>0$ such that, for any $q \in \NN$ with $1 \leq q \leq (\log N)^{\beta}$, $a \in A_q$, and $\xi = a/q + \theta \in \RR^\Gamma$, we have
\begin{align*}
\bigg|\sum_{(n,p)\in\ZZ^{k'}\times(\pm\PP)^{k''}}\ex\big(\xi\cdot\calQ(n,p)\big)\calK(n,p)\Big(\prod_{i=1}^{k''}\log|p_i|\Big)-G(a/q)\int_{\Omega}\ex\big((\xi-a/q)\cdot\calQ(t)\big)\calK(t){\rm d }t\bigg|
\\ \lesssim \big[N^{k-1}\|\calK\|_{L^{\infty}(\Omega)} \big(1 
		+
		\sum_{\gamma \in \Gamma} |\theta_\gamma| N^{|\gamma|} \big)
 + N^k\sup_{\substack{x,y \in \Omega \\ |x-y| \leq q\sqrt{k}}} |\calK(x) - \calK(y)|\big]N\exp\big(-c\sqrt{\log N}\big). 
\end{align*}
The implied constant is independent of $N,a,q,\xi$ and the kernel $\calK$.
\end{proposition}

\subsection{Multipliers for the averaging operators}
For a function $f\colon\ZZ^\Gamma \rightarrow \CC$ with finite support, we have
\[A_t^{\calQ,k',k''}f(x) = T_{\ZZ^\Gamma}[\frkM_{t}]f(x) \quad\text{and}\quad H_t^{\calQ,k',k''}f(x) = T_{\ZZ^\Gamma}[\mathfrak{n}_{t}]f(x)\]
for the discrete Fourier multipliers
\begin{equation*}
\frakm_{t}(\xi):=
\frac{1}{\vartheta_\Omega(t)}\sum_{(n,p)\in\ZZ^{k'}\times(\pm\PP)^{k''}}\ex\big(\xi\cdot\calQ(n,p)\big)\mathds{1}_{\Omega_t}(n,p)\Big(\prod_{i=1}^{k''}\log|p_i|\Big),\quad \xi\in\TT^\Gamma,
\end{equation*}
and
\begin{equation*}
\mathfrak{n}_t(\xi):=
\sum_{(n,p)\in\ZZ^{k'}\times(\pm\PP)^{k''}}\ex(\xi\cdot\calQ(n,p))K(n,p)\mathds{1}_{\Omega_t}(n,p)\Big(\prod_{i=1}^{k''}\log|p_i|\Big),\quad \xi\in\TT^\Gamma.
\end{equation*}
Their continuous counterparts are given by
\begin{equation*}
    \Phi_t(\xi):=\frac{1}{|\Omega_t|}\int_{\Omega_t}\ex(\xi\cdot\calQ(t)){\rm d} t\quad\text{and}\quad \Psi_t(\xi):={\rm p.v.}\int_{\Omega_t}\ex(\xi\cdot\calQ(t))K(t){\rm d} t
\end{equation*}
respectively. To present a unified approach, we write $M_t^{k',k''}$, $\frkY_t$, and $\Theta_t$ to represent either $A_t^{\calQ,k',k''}$, $\frkM_t$, and $\Phi_t$ or $H_t^{\calQ,k',k''}$, $\frkN_t$, and $\Psi_t$ respectively. 
 We now present the key properties of our multiplier operators that will be used in the proof of Theorem~\ref{maintheorem2}. Let $N_n:=\lfloor 2^{n^\tau} \rfloor$ for $n\in\NN$ and some $\tau\in(0,1]$ adjusted later.

\begin{enumerate}[label={\bf Property \arabic*.}, ref=\arabic*, itemindent=*,leftmargin=0pt]
\item\label{pr:1} For each $\alpha > 0$, there is $\beta_{\alpha} > 0$ such that, for any
	$\beta > \beta_{\alpha}$ and $n \in \NN$, if there is a multi-index $\gamma_0 \in \Gamma$ with
	\[
    	\bigg|\xi_{\gamma_0} - \frac{a}{q} \bigg| \leq \frac{1}{q^2}
	\]
	for some coprime integers $a$ and $q$ with $1 \leq a \leq q$ and
	$(\log N_n)^\beta \leq q \leq N_n^{\abs{\gamma_0}} (\log N_n)^{-\beta}$, then
	\[
		|(\frkY_{N_n} - \frkY_{N_{n-1}})(\xi)| \lesssim C(\log N_n)^{-\alpha}.
	\]
 This follows from Proposition~\ref{weylinq} with $\phi(x)\equiv(\vartheta_\Omega(N_n))^{-1}$ for the $\frkY_t = \frkM_{t}$ case and with $\phi(x)=K(x)$ for the $\frkY_t = \frkN_{t}$ case, noting the size condition \eqref{eq:size-unif} and the continuity condition \eqref{eq:K-modulus-cont}.
\item\label{pr:2}
Let $A$ be the $|\Gamma| \times |\Gamma|$ diagonal matrix with
\begin{equation}
	\label{matrixA}
	(A v)_\gamma = \abs{\gamma} v_\gamma.
\end{equation}
For any $t>0$, we set $t^A v: = \big(t^{\abs{\gamma}} v_\gamma : \gamma \in \Gamma\big).$ Then
\begin{equation*}
\big|\Theta_{N_n}(\xi) - \Theta_{N_{n-1}}(\xi)\big|
		\lesssim
		\min\big\{|N_n^A\xi|_\infty, |N_n^A \xi|_\infty^{-1/|\Gamma|}\big\},\quad\text{for each }n\in\NN.
\end{equation*}
In the $\Theta_t = \Phi_t$ case, this follows from the mean value theorem and the standard van der Corput lemma. In the $\Theta_t = \Psi_t$ case, this follows from the cancellation condition \eqref{eq:cancel} and \cite[Proposition B.2]{MSZ2} (see \cite[p. 21]{MSZ2} for details).

\item\label{pr:3}
For each $\alpha > 0$,  $n \in \NN$, and $\xi \in \TT^\Gamma$ satisfying
\[
    \bigg|\xi_\gamma - \frac{a_\gamma}{q} \bigg| \leq N_n^{-\abs{\gamma}} L\qquad\text{for all }\gamma \in \Gamma
\]
with $1 \leq q \leq L$, $a\in A_q^\Gamma$, and $1 \leq L \leq \exp\big(c\sqrt{\log{N_n}}\big) (\log N_n)^{-\alpha}$, we have
\[
        \frkY_{N_n}(\xi) - \frkY_{N_{n-1}}(\xi) =
        G(a/q) 
		\big(\Theta_{N_n}(\xi - a/q) - \Theta_{N_{n-1}}(\xi - a/q)\big) 
		+ \mathcal{O}\big((\log N_n)^{-\alpha}\big),
\]
for some constant $c>0$ which is independent of $n, \xi, a$ and $q$. 

In the $\frkY_t = \frkM_t$, $\Theta_t = \Phi_t$ case, this is \cite[Property 6]{Troj}. In the $\frkY_t = \frkN_t$, $\Theta_t = \Psi_t$ case, this follows from Property~\ref{pr:1} alongside Lemma~\ref{approxlemma} with $\Omega:=\Omega_{N_n}\setminus\Omega_{N_{n-1}}$ and $\calK(n,p):=K(n,p)\ind{\Omega}$, noting the size condition \eqref{eq:size-unif} and the continuity condition \eqref{eq:K-modulus-cont}. For details see \cite[Lemmas 3 and 5]{Troj}.
\end{enumerate}

\subsection{Parameters discussion}
Let $p\in(1,\infty)$ be fixed and let $\chi\in(0,1/10)$. Fix $\tau$ with $0 < \tau < 1-\min(2,p)^{-1}$ and let $N_n:=\lfloor 2^{n^\tau} \rfloor$ for $n\in\NN$. If $p \in (1,2)$, fix $p_0$ such that $1 < p_0 < p$. If instead $p \in (2,\infty)$, fix $p_0 > p$. If $p=2$, the discussion is moot since all the interpolation arguments in the article become unnecessary. We choose $\rho$ with 
\[\rho > \frac{1}{\tau}\frac{pp_0-2p}{2p_0 - 2p}\]
so that interpolation of the estimates
\[\|T\|_{\ell^2} \lesssim n^{-\rho \tau} \quad \text{and} \quad \|T\|_{\ell^{p_0}}\lesssim 1 \]
yields 
\[\|T\|_{\ell^p} \lesssim n^{-(1+\varepsilon)} \text{ for some } \varepsilon > 0.\]
Property~\ref{pr:1} gives us a corresponding $\beta_\rho$. We fix a choice of $\beta > \beta_\rho$ and then fix a choice of $u\in\NN$ with $u>|\Gamma|\beta$. We also have the value of $\delta$ coming from the Gaussian sum estimate \eqref{gausssumest}. With these fixed, we choose the value of $\varrho$ in Theorem~\ref{IW} to be
\[\varrho:= \min\bigg(\frac{\chi}{10 u}, \frac{\delta}{8\tau}\bigg).\]
\section{Proof of Theorem~\ref{maintheorem2}}

By the monotone convergence theorem and standard density arguments, it is enough to prove that
\[
\big\|\sup_{t \in \II} |M_t^{k',k''}f_{\frkI}|\big\|_{\ell^p(\ZZ^\Gamma;\ell^2)} \lesssim_{p,k,|\Gamma|} \|f_{\frkI}\|_{\ell^p(\ZZ^\Gamma;\ell^2)}
\]
holds for every finite subset $\II\subset\RR_+$ with the implicit constant independent of the set $\II$. For any $t_0 \in \II$, we have
\[\sup_{t \in \II} |M_t^{k',k''}f_{\frkI}| \leq \sup_{t \in \II} |(M_t^{k',k''}-M_{t_0}^{k',k''})f_{\frkI}| + |M_{t_0}^{k',k''}f_{\frkI}|,\]
so,
\[
\big\|\sup_{t \in \II} |M_t^{k',k''}f_{\frkI}|\big\|_{\ell^p(\ZZ^\Gamma;\ell^2)} \lesssim \calS_{\ZZ^\Gamma}^p(M_t^{k',k''}f_{\frkI}: t \in \II) + \|M_{\inf \II}^{k',k''}f_{\frkI}\|_{\ell^p(\ZZ^\Gamma;\ell^2)}.
\]
Thus, it suffices to show that 
\[\calS_{\ZZ^\Gamma}^p(M_t^{k',k''}f_{\frkI}: t \in \II) \lesssim \|f_{\frkI}\|_{\ell^p(\ZZ^\Gamma;\ell^2)}.\]

We start by splitting (cf. \cite[Lemma 1.3]{jsw}, \cite[Lemma 8.1]{MST2}) into long and short suprema along the subexponential sequence $N_n$. Letting $\II_n \coloneqq [N_n,N_{n+1})\cap\II$, we have 
\[
\calS_{\ZZ^\Gamma}^p\big(M_t^{k',k''}f_{\frkI} : t\in\II\big)\lesssim \calS_{\ZZ^\Gamma}^p(T_{\ZZ^\Gamma}[\frkY_{N_n}]f_{\frkI} : n \in \NN_0) + 
\bigg\| \Big(\sum_{n\in\NN_0} \sup_{t \in \II_n} \big|\big(M_t^{k',k''}-M_{\inf \II_n}^{k',k''}\big)f_{\frkI}\big|^2 \Big)^{1/2} \bigg\|_{\ell^p(\ZZ^\Gamma;\ell^2)}.
\]
\subsection{Short suprema}
Let $s_{n,0} < s_{n,1} < \ldots < s_{s,J(n)}$ be the increasing enumeration of $[N_n,N_{n+1}] \cap \II$ and let $r = \min(2,p)$.  Monotonicity of $\ell^p$ norms, Minkowski's inequality, and the triangle inequality give
\begin{equation}
\label{eq:short}
\bigg\| \Big(\sum_{n\in\NN_0} \sup_{t \in \II_n} \big|\big(M_t^{k',k''}-M_{\inf \II_n}^{k',k''}\big)f_{\frkI}\big|^2 \Big)^{1/2} \bigg\|_{\ell^p(\ZZ^\Gamma;\ell^2)}
\leq  \bigg(\sum_{n\in\NN_0} \bigg(\sum_{j=1}^{J(n)}\Big\| \big(M_{s_{n,j}}^{k',k''} - M_{s_{n,j-1}}^{k',k''}\big) f_{\frkI}\Big\|_{\ell^p(\ZZ^\Gamma;\ell^2)}\bigg)^r \bigg)^{1/r}.
\end{equation}
Since
\[\Big\| \big(M_{s_{n,j}}^{k',k''} - M_{s_{n,j-1}}^{k',k''}\big) f\Big\|_{\ell^p(\ZZ^\Gamma)} \leq \big\| \check{\frkY}_{s_{n,j}} - \check{\frkY}_{s_{n,j-1}}\big\|_{\ell^1(\ZZ^\Gamma)} \big\| f\big\|_{\ell^p(\ZZ^\Gamma)}\]
by Young's convolution inequality, \eqref{MZbasic} gives
\[\Big\| \big(M_{s_{n,j}}^{k',k''} - M_{s_{n,j-1}}^{k',k''}\big) f_{\frkI}\Big\|_{\ell^p(\ZZ^\Gamma;\ell^2)} \leq \big\| \check{\frkY}_{s_{n,j}} - \check{\frkY}_{s_{n,j-1}}\big\|_{\ell^1(\ZZ^\Gamma)} \big\| f_{\frkI}\big\|_{\ell^p(\ZZ^\Gamma;\ell^2)}.\] Therefore, we control the right hand side of \eqref{eq:short} by 
\begin{align*}
    \Big(\sum_{n\in\NN_0} \bigg(\big\| \sum_{j=1}^{J(n)} |\check{\frkY}_{s_{n,j}} - \check{\frkY}_{s_{n,j-1}}|\big\|_{\ell^1(\ZZ^\Gamma)} \bigg)^r \Big)^{1/r} \big\| f_{\frkI}\big\|_{\ell^p(\ZZ^\Gamma;\ell^2)} 
 \lesssim \Big(\sum_{n\in\NN_0} (n^{-r(1-\tau)} \Big)^{1/r} \big\| f_{\frkI}\big\|_{\ell^p(\ZZ^\Gamma;\ell^2)} 
 \lesssim \big\|f_{\frkI}\big\|_{\ell^p(\ZZ^\Gamma;\ell^2)}.
\end{align*}
The last estimates follow from \cite[eq. 4.2]{MeS} with $f = \delta_0$ and the discussion thereafter.

\subsection{Long suprema and the circle method}
Let $\eta\colon\RR^\Gamma \rightarrow [0,1]$  be a smooth function with
\[
	\eta(x) = 
	\begin{cases}
		1 & \text{if } |x|_\infty \leq \tfrac{1}{32 |\Gamma|}, \\
		0 & \text{if } |x|_\infty \geq \tfrac{1}{16 |\Gamma|}. 
	\end{cases}
\]
For $N\in\RR_+$, we define the scaling notation
\[
\eta_N(\xi):= \eta\big(2^{N \cdot A- N^{\chi}\cdot \rm{Id}}\xi\big)
\]
where $A$ is the matrix given in \eqref{matrixA} and $\rm{Id}$ is the $|\Gamma| \times |\Gamma|$ identity matrix. For dyadic integers $s \in 2^{u\NN}$, we define the \textit{annuli sets of fractions} by
\begin{equation}\label{annsets}
\Sigma_s := \begin{cases} \Sigma_{\leq s} & \text{ if } s=2^u, \\ \Sigma_{\leq s} \setminus \Sigma_{\leq s/2^u} & \text{ if } s>2^u, \end{cases}
\end{equation}
where the $\Sigma_{\leq \cdot}$ are the sets of Ionescu--Wainger fractions as in Theorem~\ref{IW}.
For $t\geq 2^u$, we set $F(t):=\max\{s \in 2^{u\NN} : s \leq t\}$.
We define
\[
\Xi_{\leq j^{\tau u}}(\xi) :=\sum_{a/q \in \Sigma_{\leq F(j^{\tau u})}} \eta_{j^{\tau}}(\xi - a/q)
\] 
and, for $s \in 2^{u\NN}$, we define the \textit{annuli functions}
\begin{equation}\label{annmulti}
\Xi_j^s(\xi):= \sum_{a/q \in \Sigma_{s}} \eta_{j^{\tau}}(\xi - a/q).
\end{equation}
By \eqref{annsets}, we have the telescoping property
\[
\Xi_{\leq j^{\tau u}}=\sum_{\substack{s \in 2^{u\NN} \\ s \leq j^{\tau u}}} \Xi_j^s.
\]
Note that $\eta_{j^\tau}(\xi)$ satisfies the hypothesis about the support for $\frkM$ in Theorem~\ref{IW} since $\frac{1}{8 |\Gamma|}2^{-j^{\tau}+j^{\tau\chi}} \leq e^{-j^{2 \tau u \varrho}}$ provided that $\varrho \leq \chi/(10 u)$. Using the $\Xi_{\leq j^{\tau u}}$ functions, we bound the long suprema by
\begin{align*}
\calS_{\ZZ^\Gamma}^p\Big(\sum_{j=1}^n T_{\ZZ^\Gamma}[(\frkY_{N_j} - \frkY_{N_{j-1}})\Xi_{\leq j^{\tau u}}]f_{\frkI} : n \in \NN\Big) + \calS_{\ZZ^\Gamma}^p\Big(\sum_{j=1}^n T_{\ZZ^\Gamma}[(\frkY_{N_j} - \frkY_{N_{j-1}})(1-\Xi_{\leq j^{\tau u}})]f_{\frkI} : n \in \NN\Big). 
\end{align*}
These terms correspond to major and minor arcs respectively. 

\subsection{Minor arcs}
Monotonicity of $\ell^p$ norms, telescoping, and the triangle inequality give 
\begin{align*}
\calS_{\ZZ^\Gamma}^p\Big(\sum_{j=1}^n T_{\ZZ^\Gamma}[(\frkY_{N_j} - \frkY_{N_{j-1}})(1-\Xi_{\leq j^{\tau u}})]f_{\frkI} : n \in \NN\Big) 
\leq \sum_{n=1}^\infty \big\| T_{\ZZ^\Gamma}[(\frkY_{N_n} - \frkY_{N_{n-1}})(1-\Xi_{\leq n^{\tau u}})]f_{\frkI} \big\|_{\ell^p(\ZZ^\Gamma;\ell^2)}.
\end{align*}
It then suffices to show that
\begin{equation*}
    \big\| T_{\ZZ^\Gamma}[(\frkY_{N_n} - \frkY_{N_{n-1}})(1-\Xi_{\leq n^{\tau u}})]f_{\frkI} \big\|_{\ell^{p}(\ZZ^\Gamma;\ell^2)} \lesssim n^{-(1+\varepsilon)}\|f_{\frkI}\|_{\ell^p(\ZZ^\Gamma;\ell^2)}
\end{equation*}
for some $\varepsilon > 0$. This follows from 
\begin{equation*}
    \big\| T_{\ZZ^\Gamma}[(\frkY_{N_n} - \frkY_{N_{n-1}})(1-\Xi_{\leq n^{\tau u}})]f \big\|_{\ell^{p}(\ZZ^\Gamma)} \lesssim n^{-(1+\varepsilon)}\|f\|_{\ell^p(\ZZ^\Gamma)}
\end{equation*}
by \eqref{MZbasic}. This uses Property~\ref{pr:1} and follows from the proof of \cite[Eqs. (5.8), (5.9)]{Troj} with only small changes due to our differing scaling in the definition of $\eta_N(\xi)$. We omit the details.  
\subsection{Introduction to major arcs}
Using the annuli multipliers \eqref{annmulti} and the triangle inequality, we bound the major arcs term by
\begin{align*}
\calS_{\ZZ^\Gamma}^p\Big(\sum_{j=1}^n \sum_{\substack{s \in 2^{u\NN} \\ s\leq j^{\tau u}}} T_{\ZZ^\Gamma}[(\frkY_{N_j} - \frkY_{N_{j-1}})\Xi_j^s]f_{\frkI} : n \in \NN\Big)
 \leq \sum_{s \in 2^{u \NN}} \calS_{\ZZ^\Gamma}^p\Big(\sum_{\substack{1 \leq j \leq n \\ j \geq s^{1/(\tau u)}}}  T_{\ZZ^\Gamma}[(\frkY_{N_j} - \frkY_{N_{j-1}})\Xi_j^s]f_{\frkI} : n \geq s^{1/\tau u} \Big).
\end{align*}
It then suffices to show for large $s\in 2^{u\NN}$ that
\begin{equation}\label{eq:major1}
\calS_{\ZZ^\Gamma}^p\bigg(\sum_{\substack{1 \leq j \leq n \\ j \geq s^{1/(\tau u)}}} T_{\ZZ^\Gamma}[(\frkY_{N_j} - \frkY_{N_{j-1}})\Xi_j^s]f_{\frkI} : n \geq s^{1/\tau u}\bigg) \lesssim s^{-\varepsilon} \|f_{\frkI}\|_{\ell^p(\ZZ^\Gamma;\ell^2)}
\end{equation}
for some $\varepsilon > 0$ since $\sum_{s \in 2^{u\NN}} s^{-\varepsilon} < \infty$. Let $\kappa_s:=s^{2 \lfloor \varrho \rfloor}$. By splitting the left hand side of \eqref{eq:major1} at $n \approx 2^{\kappa_s}$ into small and large scales, it suffices to prove that
\begin{equation}\label{small}
\calS_{\ZZ^\Gamma}^p\Big(\sum_{\substack{1 \leq j \leq n \\ j \geq s^{1/(\tau u)}}} T_{\ZZ^\Gamma}[(\frkY_{N_j} - \frkY_{N_{j-1}})\Xi_j^s]f_{\frkI} : n^\tau \in [s^{1/u}, 2^{\kappa_s+1}] \Big) \lesssim s^{-\varepsilon}\|f_{\frkI}\|_{\ell^p(\ZZ^\Gamma;\ell^2)}
\end{equation}
and
\begin{equation}
\label{large}
\calS_{\ZZ^\Gamma}^p\Big(\sum_{\substack{1 \leq j \leq n \\ j \geq 2^{\kappa_s/\tau}}}  T_{\ZZ^\Gamma}[(\frkY_{N_j} - \frkY_{N_{j-1}})\Xi_j^s]f_{\frkI} : n^\tau > 2^{\kappa_s} \Big) \lesssim s^{-\varepsilon} \|f_{\frkI}\|_{\ell^p(\ZZ^\Gamma;\ell^2)}.
\end{equation}

For the small scales \eqref{small}, we will use the Rademacher--Menshov inequality \eqref{eq:remark3} and Theorem~\ref{IW}. For the large scales \eqref{large}, we will use the Magyar--Stein--Wainger sampling principle from \cite[Proposition 2.1]{MSW} and its counterpart for the jump inequality from \cite[Theorem 1.7]{MSZ1}. We first recall an approximation lemma to replace our discrete multipliers with continuous counterparts. Let
\begin{equation}
    v_j^s(\xi):=\sum_{a/q \in \Sigma_s} G(a/q)\big(\Theta_{N_j} - \Theta_{N_{j-1}}\big)(\xi - a/q)\eta_{j^\tau}(\xi-a/q)
\end{equation}
and
\begin{equation}
    \Lambda_j^s(\xi):= \sum_{a/q \in \Sigma_s} \big(\Theta_{N_j} - \Theta_{N_{j-1}}\big)(\xi - a/q)\eta_{j^\tau}(\xi-a/q).
\end{equation}

\begin{lemma} \cite[Lemma 8]{MeS}
\label{multiplierapprox}
Let $M \in \NN$, $\alpha' > 0$, and $S_M:=\lfloor 2^{M^\tau - 3M^{\tau \chi}}\rfloor$. For $j\in\NN$ with $s^{1/(\tau u)} \leq j$ and $M \leq j \leq 2M$, we have
\begin{equation}
\label{approx1}
\|(\frkY_{N_j}-\frkY_{N_{j-1}})\Xi_j^s - v_j^s\|_{\ell^\infty(\TT^\Gamma)} \lesssim j^{-\alpha' \tau}
\end{equation}
and
\begin{equation}
\label{approx2}
\|(\frkY_{N_j}-\frkY_{N_{j-1}})\Xi_j^s - \Lambda_j^s \frkM_{S_M}\|_{\ell^\infty(\TT^\Gamma)} \lesssim j^{-\alpha' \tau}.
\end{equation}
\end{lemma}
\subsection{Small scales}
 Splitting $[s^{1/u},2^{\kappa_s+1}]$ into dyadic intervals and preparing via the triangle inequality to use \eqref{approx2}, we bound the left hand side of \eqref{small} by
\begin{align*}
\overbrace{\sum_{M \in 2^{\NN} \cap [s^{1/u},2^{\kappa_s}]} \calS_{\ZZ^\Gamma}^p\Big(\sum_{\substack{1 \leq j \leq n \\ j \geq s^{1/(\tau u)}}} T_{\ZZ^\Gamma}[\Lambda_j^s \frkM_{S_M}]f_{\frkI} : n^\tau \in [M,2M] \Big)}^{\text{Main Term 1}}
\\ +\underbrace{\sum_{M \in 2^{\NN} \cap [s^{1/u},2^{\kappa_s}]} \calS_{\ZZ^\Gamma}^p\Big(\sum_{\substack{1 \leq j \leq n \\ j \geq s^{1/(\tau u)}}} T_{\ZZ^\Gamma}[(\frkY_{N_j} - \frkY_{N_{j-1}})\Xi_j^s - \Lambda_j^s \frkM_{S_M}]f_{\frkI} : n^\tau \in [M,2M] \Big)}_{\text{Error Term 1}}.
\end{align*}
For Error Term~1, it will suffice to show that 
\begin{equation*}\label{error1}
    \big\|T_{\ZZ^\Gamma}[(\frkY_{N_n} - \frkY_{N_{n-1}})\Xi_n^s - \Lambda_n^s \frkM_{S_M}]f_{\frkI}\big\|_{\ell^p(\ZZ^\Gamma;\ell^2)} \lesssim n^{-(1+\varepsilon')} \|f_{\frkI}\|_{\ell^p(\ZZ^\Gamma;\ell^2)}
\end{equation*}
for some $\varepsilon' > 0$ since we would then bound it by
\begin{align*}
\sum_{n \geq s^{1/(\tau u)}} n^{-(1+\varepsilon')} \|f_{\frkI}\|_{\ell^p(\ZZ^\Gamma;\ell^2)}
\lesssim s^{-\varepsilon'/(\tau u)} \|f_{\frkI}\|_{\ell^p(\ZZ^\Gamma;\ell^2)} \lesssim s^{-\varepsilon}\|f_{\frkI}\|_{\ell^p(\ZZ^\Gamma;\ell^2)}.  
\end{align*}
This follows from 
\begin{equation*}
    \big\|T_{\ZZ^\Gamma}[(\frkY_{N_n} - \frkY_{N_{n-1}})\Xi_n^s - \Lambda_n^s \frkM_{S_M}]f\big\|_{\ell^p(\ZZ^\Gamma)} \lesssim n^{-(1+\varepsilon')} \|f\|_{\ell^p(\ZZ^\Gamma)}
\end{equation*}
by \eqref{MZbasic}, and that is \cite[Eq. 4.12]{MeS}.
 
For Main Term~1, we apply the Rademacher--Menshov inequality \eqref{eq:remark3} to bound it by
\begin{align*}
 \sum_{M \in 2^{\NN} \cap [s^{1/u},2^{\kappa_s}]} \sum_{i = 0}^{\log_2(2M)} \bigg\|\Big(\sum_j \Big| \sum_{k \in I_{i,j}^M} T_{\ZZ^\Gamma}[\Lambda_k^s \frkM_{S_M}]f_{\frkI} \Big|^2\Big)^{1/2}\bigg\|_{\ell^p(\ZZ^\Gamma;\ell^2)},
\end{align*}
where $j$ is taken over $j \geq 0$ such that $I_{i,j}^M:= [j2^i,(j+1)2^i] \cap [M^{1/\tau}, (2M)^{1/\tau}] \neq \emptyset$. Let $\tilde{\eta}_{N}(\xi):=\eta_N(\xi/2)$. Then $\tilde{\eta}_{N} \eta_{k^\tau} = \eta_{k^\tau}$ for $k^\tau \geq N$ due to the nesting supports. This lets us write
\[\Lambda_k^s \frkM_{S_M} = \Lambda_k^s \frkM_{S_M} \sum_{a/q \in \Sigma_s} \tilde{\eta}_M(\xi - a/q)=:\Lambda_k^s \frkM_{S_M}\tilde{\Xi}_{M^{1/\tau}}^s\]
for $k \in I_{i,j}^M$ since then $k \geq M^{1/\tau}$. 

By \eqref{MZgen}, it suffices to get an appropriate estimate for 
\[\bigg\|\sum_j \sum_{k \in I_{i,j}^M} \epsilon_j(\omega) T_{\ZZ^\Gamma}[\Lambda_k^s \frkM_{S_M}]f \bigg\|_{\ell^p(\ZZ^\Gamma)}\]
for any Rademacher sequence $\epsilon = (\epsilon_j(\omega))$ with $\epsilon_j(\omega) \in \{-1,1\}$ and for every $\omega \in [0,1]$.

We get the appropriate bound on $\ell^p(\ZZ^\Gamma)$ by the Ionescu-Wainger theorem and the bound for the continuous analogue 

\[\bigg\|\sum_j \sum_{k \in I_{i,j}^M} \epsilon_j(\omega) T_{\ZZ^\Gamma}[(\Theta_{N_j}-\Theta_{N_{j-1}})\eta_j^\tau]g\bigg\|_{L^p(\RR^\Gamma)} \lesssim \|g\|_{L^p(\RR^\Gamma)}\]
with a bound independent of the Rademacher sequence $\epsilon$, see \cite[Chapter XI]{bigs} or \cite{DF}. Therefore, 

\begin{equation}\label{RM_term_1}
   \bigg\|\sum_j \sum_{k \in I_{i,j}^M} \epsilon_j(\omega) T_{\ZZ^\Gamma}[\Lambda_k^s \frkM_{S_M}]f \bigg\|_{\ell^{p_0}(\ZZ^\Gamma)} 
\lesssim \big\|T_{\ZZ^\Gamma}[\frkM_{S_M}]f\big\|_{\ell^{p_0}(\ZZ^\Gamma)} \lesssim \|f\|_{\ell^{p_0}(\ZZ^\Gamma)}
\end{equation}
using the uniform $\ell^{p}$-boundedness of the averaging operators.

We get an improved bound on $\ell^2$. To do this, we use that
\[\big\|\frkM_{S_M} \tilde{\Xi}_{M^{1/\tau}}^s\big\|_{\ell^{\infty}(\TT^\Gamma)}\lesssim s^{-\delta}\]
for $M \in 2^\NN \cap [s^{1/u},2^{\kappa_s}]$, see \cite[Section 4.5]{MeS}. Then \begin{equation}
\label{RM_term_2}
    \bigg\|\sum_j \sum_{k \in I_{i,j}^M} \epsilon_j(\omega) T_{\ZZ^\Gamma}[\Lambda_k^s \frkM_{S_M}]f \bigg\|_{\ell^{2}(\ZZ^\Gamma)}
    \lesssim \big\|T_{\ZZ^\Gamma}[\frkM_{S_M}\tilde{\Xi}_{M_{1/\tau}}^s]f\big\|_{\ell^{2}(\ZZ^\Gamma)} \lesssim s^{-\delta}\|f\|_{\ell^{2}(\ZZ^\Gamma)}.
\end{equation}
Interpolation of \eqref{RM_term_1} with \eqref{RM_term_2} then gives that 
\[ 
\bigg\|\sum_j \sum_{k \in I_{i,j}^M} \epsilon_j(\omega) T_{\ZZ^\Gamma}[\Lambda_k^s \frkM_{S_M}]f \bigg\|_{\ell^{p}(\ZZ^\Gamma)}
\lesssim s^{-8\varrho}\|f\|_{\ell^p(\ZZ^\Gamma)}
\]
since $8\varrho \leq \delta/(\rho \tau)$. We then apply \eqref{MZgen}

\[\bigg\|\Big(\sum_j \Big| \sum_{k \in I_{i,j}^M} T_{\ZZ^\Gamma}[\Lambda_k^s \frkM_{S_M}]f_{\frkI} \Big|^2\Big)^{1/2}\bigg\|_{\ell^p(\ZZ^\Gamma;\ell^2)} \lesssim \|f_{\frkI}\|_{\ell^p(\ZZ^\Gamma;\ell^2)}\]

Thus, we may dominate Main Term~1 by
\begin{align*}
\sum_{M \in 2^{\NN} \cap [s^{1/u},2^{\kappa_s}]} \sum_{i = 0}^{\log_2(2M)} s^{-8\varrho}\|f_{\frkI}\|_{\ell^p(\ZZ^\Gamma;\ell^2)}\lesssim \kappa_s^2 s^{-8\varrho} \|f_{\frkI}\|_{\ell^p(\ZZ^\Gamma;\ell^2)} \lesssim s^{-4\varrho} \|f_{\frkI}\|_{\ell^p(\ZZ^\Gamma;\ell^2)}
\end{align*}
since $\kappa_s \leq s^{2\varrho}$, concluding the proof of \eqref{small}.

\subsection{Large scales} 
We bound the left hand side of \eqref{large} by
\begin{align*}
\overbrace{\calS_{\ZZ^\Gamma}^p\Big(\sum_{\substack{1 \leq j \leq n \\ j \geq 2^{\kappa_s/\tau}}}  T_{\ZZ^\Gamma}[v_j^s]f_{\frkI} : n^\tau > 2^{\kappa_s} \Big)}^{\text{Main Term 2}}
+ \overbrace{\sum_{n \geq 2^{\kappa_s/\tau}}  \big\|T_{\ZZ^\Gamma}[(\frkY_{N_n} - \frkY_{N_{n-1}})\Xi_n^s - v_n^s]f_{\frkI} \big\|_{\ell^p(\ZZ^\Gamma;\ell^2)}}^{\text{Error Term 2}}.
\end{align*}
For Error Term~2, it will suffice to show that
\begin{equation*}\label{error2}
    \big\|T_{\ZZ^\Gamma}[(\frkY_{N_n} - \frkY_{N_{n-1}})\Xi_n^s - v_n^s]f_{\frkI} \big\|_{\ell^p(\ZZ^\Gamma;\ell^2)} \lesssim e^{(|\Gamma|+1)s^\varrho} n ^{-(1+\varepsilon')} \|f_{\frkI}\|_{\ell^p(\ZZ^\Gamma;\ell^2)}
\end{equation*}
for some $\varepsilon' > 0$ since we would then bound it by
\begin{align*}
e^{(|\Gamma|+1)s^\varrho} \sum_{n \geq 2^{\kappa_s/\tau}} n^{-(1+\varepsilon')} \|f_{\frkI}\|_{\ell^p(\ZZ^\Gamma;\ell^2)}
\lesssim e^{(|\Gamma|+1)s^\varrho} 2^{-s^{2\varrho}\varepsilon'/\tau} \|f_{\frkI}\|_{\ell^p(\ZZ^\Gamma;\ell^2)} \lesssim s^{-\varepsilon} \|f_{\frkI}\|_{\ell^p(\ZZ^\Gamma;\ell^2)}.
\end{align*}
This follows from 
\begin{equation*}
    \big\|T_{\ZZ^\Gamma}[(\frkY_{N_n} - \frkY_{N_{n-1}})\Xi_n^s - v_n^s]f \big\|_{\ell^p(\ZZ^\Gamma)} \lesssim e^{(|\Gamma|+1)s^\varrho} n ^{-(1+\varepsilon')} \|f\|_{\ell^p(\ZZ^\Gamma)},
\end{equation*}
by \eqref{MZbasic}, and that is \cite[Eq. 4.15]{MeS}.

For Main Term~2, we define
\[
w^s(\xi):= \sum_{a/q \in \Sigma_s} G(a/q)\tilde{\eta}_{2^{\kappa_s}}(\xi - a/q), \quad \Pi^s(\xi):=\sum_{a/q \in \Sigma_s} \tilde{\eta}_{2^{\kappa_s}}(\xi - a/q),\]
and
\[
\omega_n^s(\xi):= \sum_{2^{\kappa_s/\tau} \leq j \leq n} (\Theta_{N_j} - \Theta_{N_{j-1}})(\xi)\eta_{j^\tau}(\xi). 
\]

Let $Q_s:=\lcm(q: a/q \in \Sigma_s).$ By property \textit{(iv)} from Theorem~\ref{IW}, we have $Q_s \leq 3^s$. The function $\omega_n^s$ is supported on $[-\frac{1}{4Q_s}, \frac{1}{4Q_s}]$ for large $s\in2^{u\NN}$ since, on the support of $\eta_{2^{\kappa_s}}$, we have $|\xi_\gamma| \leq2^{-2^{-\kappa_s}+2^{\kappa_s\chi}} \leq (4Q_s)^{-1}$
for all $\gamma \in \Gamma$ and large $s$. We also have 
\[\sum_{2^{\kappa_s/\tau} \leq j \leq n} v_j^s(\xi) = w^s(\xi)\sum_{b \in \ZZ^\Gamma} \omega_n^s(\xi - b/Q_s)  .\]
Therefore, it suffices to prove
\begin{equation}\label{largescale1}
\calS_{\ZZ^\Gamma}^p\Big(T_{\ZZ^\Gamma}\Big[\sum_{b \in \ZZ^\Gamma} \omega_n^s(\cdot - b/Q_s) \Big]f_{\frkI} : n^\tau > 2^{\kappa_s} \Big)\lesssim\|f_{\frkI}\|_{\ell^p(\ZZ^\Gamma;\ell^2)}
\end{equation}
and
\begin{equation}
\label{largescale2}
\big\| T_{\ZZ^\Gamma}[w^s]f_{\frkI}\|_{\ell^p(\ZZ^\Gamma;\ell^2)} \lesssim s^{-\varepsilon} \|f_{\frkI}\|_{\ell^p(\ZZ^\Gamma;\ell^2)}
\end{equation}
for some $\varepsilon > 0$. \eqref{largescale2} follows from
\begin{equation*}
\big\| T_{\ZZ^\Gamma}[w^s]f\|_{\ell^p(\ZZ^\Gamma)} \lesssim s^{-\varepsilon} \|f\|_{\ell^p(\ZZ^\Gamma)}
\end{equation*}
by \eqref{MZbasic}, and that is \cite[Eq. 4.17]{MeS}.

By the Magyar--Stein--Wainger sampling principle \cite[Proposition 2.1]{MSW} for the supremum seminorm, \eqref{largescale1} follows from
\begin{equation}\label{largescale1dis}
\calS_{\RR^\Gamma}^p(T_{\RR^\Gamma}[\omega_n^s]f_{\frkI} : n^\tau > 2^{\kappa_s}) \lesssim \|f_{\frkI}\|_{L^p(\RR^\Gamma;\ell^2)}.
\end{equation}
To prove \eqref{largescale1dis}, we use that the $\omega_n^s$ functions are almost telescoping. We define
\[\Delta_n^s(\xi) := \sum_{2^{\kappa_s/\tau} \leq j \leq n} (\Theta_{N_j} - \Theta_{N_{j-1}})(\xi) = (\Theta_{N_n} - \Theta_{N_{2^{\kappa_s/\tau} - 1}})(\xi). \]
Then \eqref{largescale1dis} follows from 
\begin{equation}\label{largecontdelta}
\calS_{\RR^\Gamma}^p(T_{\RR^\Gamma}[\Delta_n^s]f_{\frkI} : n^\tau > 2^{\kappa_s})\lesssim\|f_{\frkI}\|_{L^p(\RR^\Gamma;\ell^2)}
\end{equation}
since the error term is bounded by
\begin{equation*}
\sum_{n > 2^{\kappa_s/\tau}} \big\|T_{\RR^\Gamma}[(\Theta_{N_n} - \Theta_{N_{n-1}})(\eta_{n^\tau}-1)]f_{\frkI}\big\|_{L^p(\RR^\Gamma;\ell^2)}\lesssim\|f_{\frkI}\|_{L^p(\RR^\Gamma;\ell^2)}.
\end{equation*}
This last estimate follows from 
\begin{equation*}
\sum_{n > 2^{\kappa_s/\tau}} \big\|T_{\RR^\Gamma}[(\Theta_{N_n} - \Theta_{N_{n-1}})(\eta_{n^\tau}-1)]f\big\|_{L^p(\RR^\Gamma)}\lesssim\|f\|_{L^p(\RR^\Gamma)}
\end{equation*}
by \eqref{MZbasic}, and this follows from Property \ref{pr:2} and interpolation. ~\eqref{largecontdelta} itself follows from
\begin{align*}
\calS_{\RR^\Gamma}^p(T_{\RR^\Gamma}[\Theta_t]f : t > 0) \lesssim \|f\|_{L^p(\RR^\Gamma)},
\end{align*}
and this follows from \cite[Appendix A]{MST1}. 

\end{document}